\documentclass{article}

\usepackage{amssymb}
\begin{document}

\newcommand{\frg}{\mathfrak{ g}}
\newcommand{\frh}{\mathfrak{ h}}
\newcommand{\frk}{\mathfrak{k}}
\newcommand{\frp}{\mathfrak{p}}
\newcommand{\frs}{\mathfrak{s}}
\newcommand{\fru}{\mathfrak{u}}
\newcommand{\frl}{\mathfrak{l}}
\newcommand{\fro}{\mathfrak{o}}
\newcommand{\Tr}{{\rm Tr}}
\newcommand{\cA}{{\cal A}}
\newcommand{\cZ}{{\cal Z}}
\newcommand{\cB}{{\cal B}}
\newcommand{\cG}{{\cal G}}
\newcommand{\bA}{{\bf A}}
\newcommand{\Spin}{{\rm Spin}}
\newcommand{\cM}{{\cal M}}
\newcommand{\tcM}{\tilde{{\cal M}}}
\newcommand{\Sheaf}{{\cal S}}
\newcommand{\bR}{{\bf R}}
\newcommand{\bC}{{\bf C}}
\newcommand{\bP}{{\bf P}}
\newcommand{\bH}{{\bf H}}
\newcommand{\bZ}{{\bf Z}}
\newcommand{\bF}{{\bf F}}
\newcommand{\cP}{{\bf C P}}
\newcommand{\upsi}{\underline{\psi}}
\newcommand{\opsi}{\overline{\psi}}
\newcommand{\cE}{{\cal E}}
\newcommand{\ad}{{\rm ad}\ }
\newcommand{\uz}{{\underline{z}}}
\newcommand{\tA}{\tilde{A}}
\newcommand{\ccM}{\underline{M}}
\newcommand{\bE}{{\bf E}}
\newcommand{\tE}{{\tilde E}}
\newcommand{\tpi}{\tilde{\pi}}
\newcommand{\tcZ}{\tilde{{\cal Z}}}
\newcommand{\ux}{\underline{x}}
\newcommand{\uy}{\underline{y}}
\newcommand{\trho}{\tilde{\rho}}
\newcommand{\Mon}{{\rm Mon}}
\newcommand{\uMon}{\underline{{\rm Mon}}}
\newcommand{\cX}{{\cal X}}
\newcommand{\cL}{{\cal L}}
\newtheorem{prop}{{Proposition}}
\newtheorem{cor}{{ Corollary}}
\newtheorem{defn}{{Definition}}
\newtheorem{thm}{{Theorem}}
\newtheorem{formal}{{Tentative prediction}}

\title{Gauge theory in higher dimensions, II}
\author{Simon Donaldson and Ed Segal\\  Imperial College, London}
\maketitle
     

\tableofcontents

\section{Introduction}
In this paper we follow up  some of the ideas discussed in \cite{kn:DT}.
The  theme of that article was the possibility of extending  familiar constructions in gauge theory, associated to problems in low-dimensional topology, to higher
dimensional situations, in the presence of an appropriate special geometric
structure. The starting point for this was the  \lq\lq
holomorphic Casson invariant'', counting holomorphic bundles over a Calabi-Yau
3-fold, analogous to the Casson invariant which counts flat connections over
a differentiable 3-manifold. This was developed rigorously 
by Richard Thomas \cite{kn:RT}  in an algebro-geometric framework, and the ideas have been taken up by various authors \cite{kn:MNOP}.
From a differential-geometric standpoint one can make parallel discussions
of two cases: one involving gauge theory and one involving submanifolds.
There has been a considerable amount of work in a similar direction, involving
ideas of Topological Quantum Field Theory 
   \cite{kn:NCL},  \cite{kn:Salur}.

In the familiar gauge theory picture one views the Casson invariant as the
Euler characteristic of the instanton Floer homology groups. Thus it is natural
to hope for some analogous structure associated to a Calabi-Yau 3-fold. This
was discussed in a  general way in \cite{kn:DT} but  the discussion there
did not pin
down exactly what structure one could expect. That is the goal of the present
paper. In brief, we will argue that one should hope to find a {\it holomorphic
bundle} over the moduli space of Calabi-Yau 3-folds, of rank equal to the
holomorphic Casson  invariant (sometimes called the {\it DT invariant}) defined by Thomas.

Just as in \cite{kn:DT}, many of the arguments here are    tentative and speculative, since the fundamental analytical results that one
would need to develop a theory properly are not yet in place. 
These have to do with 
 the {\it compactness} of moduli spaces of solutions. While
considerable progress has been made in this direction 
 by Tian \cite{kn:T1},\cite{kn:T2}, a detailed theory---in either the
gauge theory or submanifold setting---seems still to
be fairly distant. The issues are similar to those involved in \lq\lq counting''
special Lagrangian submanifolds in Calabi-Yau manifolds, which have been considered by Joyce\cite{kn:Joyce2} but where, again, a final theory is still lacking. 

The core of this article is Section 4, where we explain how to construct
holomorphic bundles over Calabi-Yau moduli spaces, assuming favourable properties
of a \lq\lq $(6+1)$-dimensional'' differential-geometric theory. The preceding
sections 2 and 3 develop background material, mostly fairly standard but
introducing a point of view involving \lq\lq taming forms''. In Section 5
we explain how our construction matches up with standard algebraic topology,
following the familiar Floer-theory philosophy. In Section 6 we go back to
discuss  the central, unresolved,  compactness issues. We explain the relevance of recent work of Haydys which brings in a version of the \lq\lq Fueter equation''.
This perhaps points the way to a unification of the gauge theory and calibrated geometry discussions and connections with the more  algebro-geometric approach..

We are very grateful to Richard Thomas and Dominic Joyce for many discussions of this material. The paper has been substantially revised following comments of Joyce on an earlier draft. We are also grateful to Andriy Haydys for allowing us to present part of his forthcoming work. 

\section{Basic set-up}
\subsection{Exceptional holonomy}
We will begin by recalling standard material on exceptional holonomy. Some
references
are \cite{kn:Joyce1}, \cite{kn:Sa}. Start with the positive spin representation of $\Spin(7)$ on the $8$-dimensional real vector space $S^{+}$. The basic fact is that this action maps onto  the orthogonal group $SO(S^{+})=SO(8)$. Likewise for the negative spin representation. This is the phenomenon of \lq\lq triality'': there are automorphisms of $\Spin(8)$ permuting the  three representations
$S^{+}, S^{-}, \bR^{8}$.  In particular the stabiliser in $\Spin(8)$ of a unit spinor in $S^{+}$ is a copy of $\Spin(7)\subset \Spin(8)$, which maps to $\Spin(7)\subset SO(8)$.  A Riemannian $8$-manifold $X$ with a covariant constant unit spinor field has holonomy contained in $\Spin(7)$. In this situation we have a decomposition of the $2$-forms
\begin{equation}  \Lambda^{2}= \Lambda^{2}_{21}\oplus \Lambda^{2}_{7}, \end{equation}
where $\Lambda^{2}_{21}$ corresponds to the Lie algebra of $\Spin(7)$, under the isomorphism
$\Lambda^{2}= \frs\fro(8)$, and $\Lambda^{2}_{7}$ is the orthogonal complement. There is also a parallel $4$-form $\Omega$ which is equal to $(\sum \theta_{i}^{2})/7$ for any orthonormal basis $\theta_{i}$ of $\Lambda^{2}_{7}$. We can see this form in a useful explicit model.
 Suppose we have two copies
$\bR^{4}_{1}, \bR^{4}_{2}$ of $\bR^{4}$, each with spin structures. Then the positive spin space of
$\bR^{4}_{1}\oplus \bR^{4}_{2}$ is the real part of
\begin{equation}  \left(S^{+}(\bR^{4}_{1})\otimes S^{+}(\bR^{4}_{2})\right)\ \oplus\ \left( S^{-}(\bR^{4}_{1})\otimes S^{-}(\bR^{4}_{2})\right) . \end{equation}
(Recall that the spin spaces in $4$-dimensions are quaternionic and the complex tensor product of two quaternionic vector spaces has a natural real structure.)
 Fix an isomorphism $\Psi$ between $S^{+}(\bR^{4}_{1})$ and $S^{+}(\bR^{4}_{2})$. We can  regard $\Psi$  as an element of the tensor product and we get a distinguished spinor in $8$ dimensions. In other words we have a subgroup $H$ of $\Spin(7)\subset SO(8)$, locally isomorphic to $SU(2)\times SU(2)\times SU(2)$, consisting of automorphims of $\bR^{8}$ which preserve the decomposition $\bR^{4}_{1}\oplus \bR^{4}_{2}$ and $\Psi$.
In this picture the $4$-form $\Omega$ corresponding to our distinguished spinor is
\begin{equation}  dx_{1} dx_{2} dx_{3} dx_{4} + dy_{1}dy_{2} dy_{3} dy_{4} +\sum_{i=1}^{3}
\omega_{i} \wedge \omega'_{i}. \end{equation}
Here $x_{i}, y_{i}$ are standard co-ordinates on the two copies of $\bR^{4}$,
$\omega_{i}$ is a standard orthonormal basis for $\Lambda^{+}(\bR^{4}_{1})$
and $\omega'_{i}$ the basis of $\Lambda^{+}(\bR^{4}_{2})$ which corresponds to this
under the isomorphism induced by $\Psi$. In fact this form $\Omega$ determines the spinor, so we could also
define $\Spin(7)\subset GL(8,\bR)$ to be the stabiliser of this $4$-form. The $GL(8,\bR)$
orbit of $\Omega$ 
is a $43$-dimensional submanifold $\cA\subset \Lambda^{4} \bR^{8}$,
which can be viewed as $GL(8,\bR)/\Spin(7)$. On any $8$-manifold we have a
copy of $\cA$ associated to each tangent space in the obvious way and a Spin(7)
structure is equivalent to  a closed $4$-form which takes values in this subset.
  
Now consider a unit spinor in $S^{+}$ and a unit vector in $\bR^{8}$. The stabiliser of the pair is the exceptional Lie group $G_{2}$, which can be regarded as a subgroup of $SO(7)$. This means that a Riemannian product $\bR\times Y$ has holonomy contained in $\Spin(7)$ if and only if the holonomy of $Y$ is contained in $G_{2}\subset SO(7)$. On such a manifold $Y$ we have a decomposition of the $2$-forms 
\begin{equation}  \Lambda^{2}= \Lambda^{2}_{14}\oplus \Lambda^{2}_{7}\end{equation}
where $\Lambda^{2}_{14}$ corresponds to the Lie algebra of $G_{2}$. We have a covariant constant $3$-form $\phi$ and $4$-form $\sigma$, such that on the cylinder
\begin{equation}   \Omega= \phi dt + \sigma. \end{equation}

We can  define the concept of a non-compact Riemannian $8$ manifold with holonomy in $\Spin(7)$ and a  {\it tubular end} modelled on $ (0,\infty)\times
Y$, for a $G_{2}$-manifold $Y$. (That is, the metric differs from the product by an exponentially decaying term.) We can also consider a \lq\lq neck-stretching'' sequence of
$Spin(7)$-structures on a compact manifold that degenerate to a limit which is the disjoint union of two such non-compact manifolds. That is, the metrics contain regions which are almost isometric to long finite tubes
$ (-T_{i}, T_{i})\times Y$ , where $T_{i}\rightarrow \infty$. 

Now we can repeat the discussion, starting in $7$ dimensions. Considering
$G_{2}\subset SO(7)$, the stabiliser of a unit vector is a copy of $SU(3)\subset SO(6)$. A Riemannian product $\bR\times Z$ has a $G_{2}$-structure if and only if $Z$ is a {\it Calabi-Yau 3-fold}, with holonomy in $SU(3)$. Then we have a decomposition of the $2$-forms
\begin{equation}\Lambda^{2}= \Lambda^{2}_{8}\oplus \Lambda^{2}_{7}, \end{equation}
where $\Lambda^{2}_{8}$ corresponds to the Lie algebra of $SU(3)$. There are covariant constant $3$-forms $\rho_{1},\rho_{2}$ and a $2$-form $\omega$ such that on the cylinder
\begin{equation} \sigma= \rho_{2} ds + \omega^{2}\ \ , \  \phi= \omega ds + \rho_{1}. \end{equation}
Our notation here is that  $s$ is the co-ordinate on $\bR$. 
In fact either one of $\rho_{1},\rho_{2}$ determines the other since
$$  \rho_{2}= -I \rho_{1}, $$
where $I$ is the parallel complex structure. From another point of view, the complex combination $\theta=\rho_{1}+ i \rho_{2}$ is a holomorphic $3$-form on $Z$.  The $2$-form $\omega$ lies in $\Lambda^{2}_{7}$, so we get a
further decomposition
\begin{equation}  \Lambda_{2}=\Lambda^{2}_{8} \oplus \Lambda^{2}_{6}\oplus \langle \omega\rangle. \end{equation} 
Again we may consider $G_{2}$ manifolds with tubular ends and neck-stretching sequences. 

Here we stop, although we could repeat the process to consider Calabi-Yau $3$-folds with tubular ends, etc. 
This kind of neck-stretching sequence, and the converse \lq\lq gluing theory''
for manifolds with tubular ends, is central in the work of Kovalev \cite{kn:Kov},
and many interesting new examples have been found recently by Kovalev and
Nordstrom \cite{kn:KN}.

\subsection{Gauge theory and submanifolds}

Next we review slightly less standard material on auxiliary differential geometric objects: submanifolds and connections. A fundamental reference for the first is  the work of Harvey and Lawson
\cite{kn:HL}; a number of references for the second can be found in \cite{kn:DT}.

Start again in dimension $8$. Take our model $\bR_{1}^{4}\oplus
\bR^{4}_{2}$ above and consider the $\Spin(7)$-orbit of the $4$-plane $\bR^{4}_{1}$
in the Grassmannian of oriented $4$-planes in $\bR^{8}$. This is the set
of {\it Cayley $4$-planes}, and forms a 12-dimensional submanifold  in the full
Grassmannian (since the stabiliser of $\bR^{4}_{1}$ in the $21$-dimensional
group $\Spin(7)$ is  the $9$ dimensional subgroup $H$).  Another definition is that an oriented $4$-plane $\Pi$ is Cayley if the restriction of $\Omega$
to $\Pi$ is the volume form. A third is that for any vector $v$ in $\bR^{8}$
\begin{equation}  i_{v}(\Omega)\vert_{\Pi} = *_{\Pi}( v^{\sharp}\vert_{\Pi}), \end{equation}
where $*_{\Pi}$ is the $*$ operator on $\Pi$ induced by the metric and $v^{\sharp}\in
(\bR^{8})^{*}$
is the dual of $v$, again defined by the metric. 

 Now in an $8$-manifold with a $\Spin(7)$ structure we may consider {\it Cayley submanifolds}, whose tangent space at each point is Cayley. There are two fundamental properties of this  condition:
\begin{itemize}
\item {\bf Property A} The condition is an {\it elliptic} PDE. As a check on this, note that the condition is given locally by $4$ equations, since $V$ has codimension $4$ in the full Grassmannian, while $4$-dimensional submanifolds near to a given one can be represented by sections of the four-dimensional normal bundle, so the condition is determined. In fact the linearisation of the condition is given by a version of the Dirac operator acting on sections of the normal bundle \cite{kn:McL}.
\item {\bf Property B} The submanifold is {\it calibrated}: its volume is the topological invariant given by the integral of the closed form $\Omega$ and it is a minimal
submanifold, minimizing volume in its homology class.
\end{itemize}

Next we move to gauge theory. We define a $\Spin(7)$-instanton to be a connection on a bundle $E$, with structure group $SU(l)$ say,  whose curvature lies in $\Lambda_{2}^{21}$.
Then we have, likewise:
\begin{itemize} \item {\bf Property A} The equation is elliptic, modulo gauge equivalence. The linearised theory can be formulated in terms of a bundle-valued version of the {\it elliptic complex}
\begin{equation} \Omega^{0}\stackrel{d}{\rightarrow}\Omega^{1}
\stackrel{\pi_{7}\circ d}{\rightarrow} \Omega^{2}_{7}. \end{equation}
(Again, note the dimension check: $1-8+7=0$.)
\item {\bf Property B} The  Yang-Mills energy is determined by the topology
of the bundle. This comes from the algebraic fact that for $\alpha$ in $\Lambda^{2}_{21}$ we have
$$  \alpha \wedge \alpha \wedge \Omega = -\vert \alpha \vert^{2} {\rm vol} $$
Thus for a $\Spin(7)$-instanton:
\begin{equation}   \int_{X^{8}} \vert F \vert^{2} = -\int_{X^{8}} \Tr(F^{2}) \wedge \Omega=   8\pi^{2} c_{2}(E)\cup\Omega \ [X]. \end{equation}
\end{itemize}
It follows that $\Spin(7)$-instantons are Yang-Mills connections.
This theory has been developed in the thesis of C.Lewis \cite{kn:CL}.

\

{\bf  Optimistically} we could hope that these good properties   would allow us to define {\it invariants}, modelled on the Gromov-Witten invariants of symplectic manifolds, in the submanifold case, and on the instanton invariants of $4$-manifolds in the gauge theory case. In the simplest situation, where the relevant index or formal dimension is zero, one would seek to define a numerical invariant by counting solutions with appropriate signs. When the index is positive one could try to evaluate natural cohomology classes on the moduli space of solutions.

There is a third  property--which one might call the Floer picture---of the equations, which becomes crucial when one
considers manifolds with tubular ends and neck-stretching sequences.
Consider a tube $I\times Y$, where $I\subset \bR$ is an open interval,
finite or infinite, and $Y$ is a $G_{2}$-manifold with a $4$-form $\sigma$. A connection $\bA$ over
$I\times Y$ can be viewed as  a $1$-parameter family $A_{t}$ of connections
over $Y$. 

{\bf Property C} There is a locally-defined function $f$, on the space of
connections over $Y$ modulo gauge equivalence, such that $\Spin(7)$ instantons
correspond to gradient curves of $f$, and
\begin{equation}  \frac{d}{dt} f(A_{t}) = - \Vert \frac{dA_{t}}{dt} \Vert^{2}=-\Vert
F(A_{t})\wedge \sigma \Vert^{2}. 
\end{equation}

We know from analogous Floer-type theories that this is the essential  property needed to control solutions over infinite tubes, and  to obtain uniform control
in neck-stretching sequences. The point is that one arrives  in a situation
where $f(A_{t})$ is well-defined and the variation of $f$ over the interval
is known, so the gradient property gives  bounds on
$$  \int_{I} \Vert \frac{dA_{t}}{dt}\Vert^{2} \ dt, \ \ \ , \ \ \ \int_{I}
\Vert F(A_{t})\wedge \sigma \Vert^{2}  \ \ dt. $$

 To  explain in a little  more detail, we define  a $1$-form on the space of connections over $Y$ by mapping a tangent vector $\delta A$ to
\begin{equation}   \int_{Y} \Tr( \delta A\  F ) \wedge \sigma. \end{equation}
This $1$-form arises, locally in the space of connections modulo gauge equivalence, as the derivative of a function which can be written schematically as 
\begin{equation} f(A)= \int_{Y} CS(A) \wedge \sigma, \end{equation}
where $CS$ denotes the Chern-Simons form. Of course this is not really well-defined and a more precise definition is this. We choose a base point $A_{0}$ and for any nearby connection $A$ we choose a connection $\bA$ on a bundle over $[0,1]\times Y$ with boundary values $A_{0}$ and $A$. Then we define
\begin{equation}  f(A)= \int_{[0,1]\times Y} \sigma \ \Tr F_{\bA}^{2}  . \end{equation}
Just as in the usual Floer theory over $3$-manifolds, this function is not globally well-defined, but the indeterminacy comes from the periods of $\sigma$ over $H_{4}(Y,\bZ)$. This indeterminacy is related the notion of an \lq\lq adapted bundle''
which we will discuss further in Section 4.  But in any case we have a well-defined
closed $1$-form on the space of connections over $Y$. The $\Spin(7)$-instanton
equation over the tube can written as

\begin{equation}  \frac{dA_{t}}{dt}= F(A_{t}) \wedge \sigma \end{equation}
which displays $A_{t}$ as an integral curve of the vector field dual to this
$1$-form. 
 Likewise in the submanifold set-up, we define a $1$-form on the space of 3-dimensional submanifolds of $Y$ by mapping a variation $v$---a vector field along a submanifold $P^{3}\subset Y^{7}$---to
\begin{equation}\int_{P} i_{v}(\sigma) . \end{equation}
This $1$-form is the derivative of a locally-defined function, determined by integrating $\sigma$ over cobordisms in $[0,1]\times Y$, and we have an
analogue of Property C above.

Now we consider a connection on a bundle over $Y^{7}$ whose pull-back to the cylinder is a $\Spin(7)$-instanton. Expressed directly over $Y$ this condition is just that
\begin{equation}   F\wedge \sigma=0, \end{equation}
and we call the solutions $G_{2}$-instantons. In the Floer picture, these
are viewed as the zeros of the $1$-form on the space of connections.

 A $3$-dimensional submanifold $P\subset Y$ is called associative if  $\bR\times P$ is Cayley. In the picture above, such a submanifold is  viewed  as a zero of the $1$-form
or critical point of the locally-defined function $f$. Associative submanifolds
 can be defined more directly by the condition that for any vector $v\in TY$ the restriction to $P$  of the contraction $i_{v}(\sigma)$ vanishes. There is a basic algebraic model for the tangent space of $Y$ at a point of an associative submanifold like that which we saw in $8$ dimensions. We consider a $3$-dimensional space $\bR^{3}$ and a $2$-dimensional complex vector space
$V$ with symmetry group $SU(2)$. Then the tensor product of $V$ with the spin space $S=S(\bR^{3})$ has a real form $\bR^{4}= (V\otimes S)_{\bR}$ and we have a natural isomorphism
$ \Lambda^{2}_{+}\bR^{4}= \bR^{3}$.
In other words, if $y_{1}, y_{2}, y_{3}$ are standard co-ordinates on $\bR^{3}$ we have a corresponding basis $\omega_{1},\omega_{2}, \omega_{3}$ for the self-dual $2$-forms on $\bR^{4}$.  Then we have a $3$-form
\begin{equation}  \phi= \sum dy_{i} \omega_{i} + dy_{1} dy_{2} dy_{3}, \end{equation}
and a $4$-form
\begin{equation} \sigma= \sum_{i,j,k\ {\rm cyclic}} \omega_{i} dy_{j} dy_{k} + dx_{1} dx_{2} dx_{3} dx_{4}. \end{equation}
This gives a standard model for the tangent space of a $G_{2}$-manifold at
a point of an associative submanifold, with $y_{i}$ the co-ordinates on the
submanifold.

\

A $4$-dimensional submanifold $Q$ of $Y$ which is  Cayley when regarded as a submanifold of $\bR\times Y$ is called {\it co-associative}. The condition can be defined more directly by saying that the restriction of the 3-form $\phi$ to $Q$ vanishes. 

\

If we prefer, we can forget the $8$-dimensional
geometry and start directly in $7$ dimensions, considering $G_{2}$ manifolds
and three kinds of differential geometric objects: 
associative submanifolds,  coassociative submanifolds and $G_{2}$-instantons. These three conditions enjoy the same crucial Properties A,B discussed
above.  The equations are elliptic when set up suitably, but this is now
less straightforward. At first sight, the $G_{2}$-instanton equation imposes $7=\dim \Lambda^{2}_{7}$ conditions on the curvature of a connection, whereas we would expect to only impose $7-1$, taking account of gauge invariance. The explanation for this is that the curvature of any connection satisfies the  Bianchi identity and the linearised problem can be formulated in terms of the elliptic complex
\begin{equation}  \Omega^{0}\stackrel{d}{\rightarrow}\Omega^{1}\stackrel{\sigma\wedge d}{\rightarrow} \Omega^{6}\stackrel{d}{\rightarrow} \Omega^{7}. \end{equation}
Then the whole theory of local deformations of solutions to the $G_{2}$-instanton equations  is closely analogous to the Taubes/Floer discussion of flat connections over $3$-manifolds, with (21) taking the place of the de Rham complex over a $3$-manifold. This theory has been developed in unpublished work of
A. Tomatis and in the thesis of Henrique Sa Earp \cite{kn:HSE}

The associative condition is determined elliptic, with the linearisation
given by a version of the Dirac operator acting on sections of the normal
bundle.
 The co-associative condition is  at first sight overdetermined, since it imposes $4$ constraints on the sections of the $3$-dimensional normal bundle of a $4$-dimensional submanifold of $Y^{7}$, but when set-up properly becomes elliptic \cite{kn:McL}. One way of doing this is to embed the discussion in that of Cayley submanifolds in the cylinder. Our third characterisation (10) of Cayley
subspaces shows that if $Q$ is any Cayley submanifold in $ \bR\times Y$ then
 the $\bR$-co-ordinate $t$ is a harmonic function on $Q$. Thus if $Q$
is compact it must be a co-associative submanifold in some \lq\lq fixed time'' slice. .

\

In sum, we would {\bf optimistically} hope first for a $7$-dimensional
theory, bearing on a compact $G_{2}$-manifold $Y$ and yielding invariants
counting associative submanifolds,  $G_{2}$-instantons and co-associative
submanifolds.
In the first two cases these could be viewed as generalising the Casson invariant
(since in each case we are counting the zeros of a closed $1$-form). But
second we could hope for a $(7+1)$-dimensional theory, assigning Floer groups
to $Y$ which 
 should play the same role vis-a-vis the $\Spin(7)$ discussion , for manifolds with tubular ends and neck-stretching limits, as the ordinary Floer theory does to $4$-manifold invariants.
Of course the second theory would be a refinement of the first, since we would view the Casson invariant as the Euler characteristic of the Floer groups.

\

Now we repeat the discussion, dropping dimension again. The  $G_{2}$-instantons on $\bR\times Z^{6}$ which are lifted up from $Z$ are connections with
\begin{equation} F\wedge \omega^{2}=0 \ \ , \ \ F\wedge \rho_{2}=0. \end{equation}
(The conditions $F\wedge \rho_{1}=0, F\wedge \rho_{2}=0$ are equivalent.)
From the point of view of complex geometry, the second condition is that
$F$ has type $(1,1)$, so the connection defines a holomorphic structure on the bundle. The other condition is the Hermitian Yang-Mills equation, and we know that the solutions correspond to \lq\lq polystable'' holomorphic bundles.
Likewise a product 
$\bR\times \Sigma$ is co-associative if and only if $\Sigma$ is a complex curve in $Z^{6}$. These geometric objects in $Z^{6}$ again have the same good properties: they are defined by elliptic equations and have topological volume/energy bounds.  We will take this discussion of these   $(6+1)$-dimensional
theories further in Section 4. One can also discuss co-associative submanifolds
in a similar framework \cite{kn:NCL}, \cite{kn:Salur}. The corresponding objects in 6 dimensions are the
special Lagrangian submanifolds, which can be described as the critical points
of a locally-defined functional \cite{kn:RT2}. The co-associative submanifolds
of a tube 
are the gradient curves of this functional.

\section{Taming forms }
 \subsection{Dimension 8}
We will now take another point of view,  beginning again
in $8$-dimensions. Our model is Gromov's notion of a symplectic form \lq\lq taming''
an almost-complex structure \cite{kn:Gr}.  Let $\Omega_{0}$ be the standard form on $\bR^{8}$. The
convex hull of the set of negative squares $-\alpha \wedge \alpha$, for $\alpha$ in $\Lambda^{2}_{27}$ is a proper cone $K$ in $\Lambda^{4}\bR^{8}$ so we have a
dual cone of $4$-forms $\Omega'$ such that $\Omega'\wedge \chi>0$ for non-zero
$\chi\in K$. This is equivalent to saying that the quadratic form
$$  \alpha \mapsto \alpha \wedge \alpha \wedge \Omega'$$
is positive definite on $\Lambda^{2}_{27}$.  One can check that the subset
of Cayley planes in $Gr_{4}(\bR^{8})\subset \Lambda^{4}\bR^{8}$ lies in the cone $K$, so such a $4$-form is
also strictly positive on Cayley subspaces.   Now suppose that $\Omega$ is any  $4$-form on an $8$-manifold $X$ which
lies in the preferred subspace $\cA$ at each point. We could call this an
\lq\lq almost $\Spin(7)$ structure''. Then Cayley submanifolds
and $\Spin(7)$ instantons are defined and the equations are elliptic, just
as before. However we lose the volume/energy identity for solutions. Suppose
however we have a closed $4$-form $\Omega'$ which, at each point,  lies in the dual cone above. Then for any $\Spin(7)$-instanton we have a slightly
weaker version of \lq\lq Property B'': 
\

{\bf Property $B'$} 
\

\begin{equation} \int_{X} \vert F \vert^{2} \leq -C \int_{X} \Tr(F^{2}) \wedge \Omega',
\end{equation}

\

where the right hand side is a topological invariant. Here the constant $C$
depends only on $\Omega, \Omega'$. Similarly, we get a volume bound for Cayley submanifolds. It seems reasonable to hope that the putative theory in the
case of genuine $\Spin(7)$ structures extends to this more general situation;
in any case we will assume this is so  for the purposes of our discussion.
We call such a pair $(\Omega,\Omega')$ a {\it tamed almost-$Spin(7)$ structure}
on $X$.
The advantage of this extension is that the notion is much more flexible.
On the hand, starting with a genuine Spin(7) structure,  it gives us scope to deform the equations, for example to achieve transversality. This is
essentially the use we will make of the idea in Section 4. On the
other hand, tamed almost $\Spin(7)$ structures should be much easier to construct
than $\Spin(7)$ structures since the condition on $\Omega'$ is an open condition
on a closed $4$-form.  

\subsection{Dimension 7}
We may consider dimension reductions of the theory, in this more general
setting. Let $V$ be an oriented $7$-dimensional vector space and $GL_{+}(V)$
the linear automorphisms of positive determinant. There is an open $GL_{+}(V)$ orbit
$P_{3}\subset \Lambda^{3} V^{*}$ of \lq\lq positive'' forms, each of which has stabiliser isomorphic to the compact group $G_{2}\subset SL(V)$. Similarly there is an open orbit of positive $4$-forms $P_{4}\subset \Lambda^{4} V^{*}$ and a $GL_{+}(V)$-equivariant diffeomorphism
$*:P_{3}\rightarrow P_{4}$. We also denote the inverse map by $*$. The choice of notation is derived from the fact that each element $\phi\in P_{3}$ defines a natural Euclidean metric $g_{\phi}$ on $V$ and $*\phi$ is the usual Hodge dual defined by this metric. Likewise for any $\sigma\in P_{4}$. 

Now consider an oriented $7$-manifold $Y$. At each point $p\in Y$ we have open subsets $P_{3,p}\subset \Lambda^{3} T^{*}Y_{p}, P_{4,p}\subset \Lambda^{4} T^{*}Y_{p}$. We define an {\it almost $G_{2}$ structure} on $Y$ to be a $4$-form $\sigma$ on $Y$ which lies in $P_{4,p}$ at each point $p$. Of course it is the  same to start with  a $3$-form $\phi$ which lies in $P_{3,p}$ at each point, and an  almost $G_{2}$ structure defines a Riemannian metric on $Y$.  If $\sigma$ is an almost $G_{2}$-structure then the form
 \begin{equation} \Omega= \sigma + *\sigma dt, \end{equation}
 yields an almost $\Spin(7)$ structure on the cylinder $\bR\times Y$. Let
 $\phi', \sigma'$ be respectively, a closed $3$-form and $4$-form on $Y$,
 so $\Omega'=\sigma'+\phi' dt$ is closed $4$-form on the cylinder. Then there
 is a certain open set of \lq\lq taming pairs '' $(\phi',\sigma')$ such that
 $\Omega'$ is a taming form. The discussion of associative submanifolds
 in $Y$ proceeds just as before. They are defined by the condition that
 the restriction of
 $i_{v} \sigma$ vanishes, for each tangent vector $\sigma$. The condition
 is a determined  elliptic equation of Dirac type  and their
 volume is controlled by the integral of $\phi'$ over $P$, a topological invariant.

 The discussion of $G_{2}$-instantons is  different in an important
 way. For general $\sigma$ the equation $\sigma \wedge F$ for a connection
 is overdetermined and we do not expect any solutions. However we can consider
 connections over $\bR\times Y$ with the property that the restriction to
 each slice $\{t\}\times Y$ lies in a fixed gauge equivalence class. This
 corresponds to a pair $(A,u)$ over $Y$ where $A$ is a connection on a bundle
 $E$ and $u$ is a section of the adjoint bundle $\frg_{E}$. Then we have
 a connection $A+u  dt$ on the lifted bundle over $\bR\times Y$. The $\Spin(7)$
 instanton equation, expressed in $7$-dimensions, becomes 
 \begin{equation}     F\wedge \sigma = * d_{A} u , \end{equation}
 where $*$ is the Hodge $*$-operator of the metric $g_{\sigma}$.
  This is
 a determined elliptic equation (modulo gauge equivalence) for the pair $(A,u)$.
 This set-up is similar to one considered in early work of Thomas: the
 advantage in our present situation is that, given a taming pair $(\sigma',
 \phi')$, we still get a topological bound like (23) on the Yang-Mills energy.
   (Although $\sigma'$ does not appear explicitly  we need to
   use it  in deriving the inequality.)

   Suppose however that we assume that $\sigma$ is closed over $Y$ (in which
   case we might prefer to restrict attention to taming pairs with $\sigma'=\sigma$).
    Then
   the Bianchi identity implies that for any solution $(A,u)$ as above (over
   a compact manifold $Y$) we
   have $d_{A}u=0$. Thus in this case we do get a good theory of $G_{2}$-instantons,
   without the extra field $u$.
   The explanation is that when $\sigma$ is closed we have a locally defined
   \lq\lq Chern-Simons functional'' just as before, with critical points
   the $G_{2}$-instantons, and the linearised theory
   can be expressed in terms of an elliptic complex (21). Moreover we have exactly
   the same  gradient curve description as before. To sum up
   
   \
   
   (1) If we restrict
   attention to tamed structures on $8$-manifolds with tubular ends such
   that on  each end the structure is defined by a closed $4$-form $\sigma$
   on the cross section than we get a $(7+1)$-dimensional differential geometric
   theory with the good
   Properties$ A,B',C$. So we have reason to hope that some kind of Floer
   theory can be introduced into this more general, flexible, situation. The
   same applies in the submanifold setting, for Cayley submanifolds with
   \lq\lq associative limits''. 
   
   (2) If on the other hand we are just interested in a compact $7$-manifold
   $Y$ we can study solutions $(A,u)$ of equation (25) for any positive $4$-form
   $\sigma$, not necessarily closed, so long as there  are taming forms.
    This equation has properties
   $A$ (with Fredholm index $0$) and $B'$, so we expect to define a numerical
   invariant. This can be regarded as counting the zeros of a vector field
   on the space of connections modulo gauge, but the vector field is not
   dual to a closed $1$-form. Similarly for aasociative submanifolds.

   \
   
  {\bf Remark} 
 We note in passing that, by contrast, we expect a good theory of {\it coassociative}
 submanifolds for almost $G_{2}$ structures where the $3$-form $*\sigma$
 is closed, rather than the $4$-form $\sigma$. The problem is that, for general
 $\sigma$ the equations are overdetermined so we lose Property A. The condition
 $d*\sigma=0$ is the integrability condition for the overdetermined system, much as the condition $d\sigma=0$ is the integrability condition for the $G_{2}$-instanton equations. 
 One way of seeing this is to embed the discussion in that of Cayley submanifolds
 in the cylinder, as before.
 
 \
 
 \subsection{Dimension 6}
 The discussion becomes considerably more confusing when we go down to $6$
 dimensions. We will need to have closed  taming forms, comprising a $2$-form
 $\omega'$ a pair of $3$-forms $\rho'_{1}, \rho'_{2}$ and a $4$-form $\tau'$.
 Thus on $\bR^{2}\times Z$ our taming form will be
 $$\Omega'= \omega' ds dt + \rho'_{1} dt + \rho'_{2} ds + \tau'. $$
 Then we will consider an $SU(3)$-structure on the tangent bundle. We recall that there is an open $GL(6,\bR)$-orbit of \lq\lq positive'' $3$-forms in $\Lambda^{3}\bR^{6}$ each of which determines an almost-complex structure. Thus we can take  our $SU(3)$-structure to be specified by  a $3$-form $\rho_{1}$, which defines an almost
 complex structure and hence $\rho_{2}=-I\rho_{1}$ and a $2$-form $\omega$
 which is a positive form of type $(1,1)$ with respect to this almost complex
 structure. The point to emphasise that, in this most general formulation,
 there are algebraic constraints on the forms $\omega,\rho_{1}$ but only
 open conditions on the taming forms. Then we get a wide variety of different
 extra conditions we can impose, intermediate between this most general formulation
 and the case of genuine Calabi-Yau structures.  
 
 \

 Start with a $7$-manifold $Y$ with tubular ends and  a tamed
 almost $G_{2}$-structure defined by a $4$-form $\sigma$, not necessarily
 closed. Then we want to study the equation (25),  for pairs $(A,u)$ on $Y$.  We need to have Property C: a gradient description on the ends. For this
 we consider a closed $3$-form $\rho$ and the closed $4$-form $\tau$
 on the cross-section $Z$ and the functional
 $$  f(a,u) =\int_{Z} CS(a)\wedge \rho + \Tr( u  F(a))\wedge \tau, $$
 on pairs $(a,u)$ over $Z$. Suppose for the moment that we take an arbitrary
 Riemannian metric on $Z$ then we have a gradient equation
\begin{equation}\frac{da}{ds}= *(F\wedge \rho + d_{a}u\wedge \tau)\ ,\
 \frac{du}{ds}=
*(F\wedge\tau)
\end{equation}
On the other hand, if we can write the $4$-form $\sigma$ as
$\sigma=\rho_{2} ds +\omega^{2}$  the equation (25) becomes
\begin{equation} \omega^{2} \wedge \frac{da}{ds}= \rho_{2} \wedge F + * d_{a}u \
\ , \ \   \frac{du}{ds}= *(F\wedge \omega^{2}) \end{equation}
So we need to arrange that these equations are the same. First we should take $\tau=\omega^{2}$, that is we should suppose that $\omega^{2}$ is closed.
Second we should suppose that $\rho=I\rho_{2}$ for the almost complex structure defined by $\rho_{2}$, and we should use the standard metric associated to $\omega$ and this almost-complex structure. Then the equations (26) and (27) are equivalent.  But we can also use $\rho$ to define the same complex structure, and hence $\rho_{2}$.
So if we start on $Z$ with a closed positive $3$-form $\rho$ and a $2$-form $\omega$ such that $\omega^{2}$ is closed and $\omega$ is positive of type $(1,1)$ with respect to the almost-complex structure defined by by $\rho$ we get a model for a tubular end on which equation (25) has a gradient description  (Property C). We can take $\rho,\tau$ as two of our taming
 forms $\rho'_{1},\tau'$ and we also need another taming $3$-form $\rho'=\rho'_{2}$ and a taming $2$-form $\omega'$. We might want to suppose that in fact $\omega'=\omega$. (The type $(1,1)$ condition is the same as saying that $\omega\wedge \rho=0$.) 
In this way we obtain a $(6+1)$-dimensional theory with Properties $A,B', C$.
The stationary solutions over tubes correspond to pairs $(a,u)$ on $Z$ with
\begin{equation} F(a)\wedge \rho= d_{a}u\wedge \omega^{2} \ \ \ F\wedge \omega^{2}=0.
\end{equation} 
The data we need over $Z$ consists of closed forms $\rho,\rho',\omega$ but
the only identity (as opposed to open condition) that  we need to impose is $\rho\wedge \omega=0$. We expect that
we then get Floer groups associated to $Z$,  related to the invariants of
almost $G_{2}$-structures in $7$-dimensions.

If we want to study the $7$-dimensional theory with $\sigma$ closed, which
fits into the $(7+1)$ dimensional discussion, we are much more restricted. Then
we need  both $\rho$ and $I\rho$ to be closed $3$-forms on $Z$, which can
only occur if we have a genuine Calabi-Yau structure. In this case an integration
by parts shows that for any solution of (28), $d_{a}u$ vanishes so we are back to the
equations (22). In the opposite direction, if we consider the most general
set up with no particular relation between the closed taming $3$-forms and
the $SU(3)$ structure, we can
consider the system of equations for a triple $(a,u,v)$ where $a$ is a connection and
$u,v$ are sections of $\frg_{E}$:
\begin{equation} F\wedge \rho= d_{a} u \wedge \omega^{2}\ \ ,\ \  F\wedge (I\rho)=
d_{a}v \wedge \omega^{2}\ \ , \ \ F\wedge \omega^{2}= [u,v]\omega^{3} \end{equation}
We get a theory with Properties $ A,B'$ and we xpect that counting solutions will
generalise the holomorphic Casson invariant to
this situation.  

\

\

 Much of the above has a good  formal interpretation.  The space of connections on a bundle $E\rightarrow
 Z$ has a symplectic form
 $$  \langle \delta_{1}, \delta_{2} \rangle= \int_{Z} \Tr(\delta_{1}\wedge \delta_{2}) \wedge \omega^{2}. $$
 The symplectic quotient $\cA//\cG$ by the gauge group is given by solutions of the equation $F\wedge \omega^{2}=0$, modulo gauge equivalence. We have
an induced symplectic form on $\cA//\cG$. Any closed $3$-form $\trho$ yields
a locally-defined function on the space $\cA/\cG$
$$  f_{\trho}(a)=\int_{Z} CS(a) \wedge \trho\ . $$
The equation $F\wedge \trho= d_{a} u \wedge \omega^{2}$ is the condition defining
a critical point of $f_{\trho}$, {\it restricted to the symplectic quotient}
$\cA//\cG\subset \cA/\cG$. The evolution equation
$$  \frac{da}{dt} \wedge \omega^{2}= F\wedge \trho, $$
defines the {\it Hamiltonian flow} associated to the function $f_{\trho}$
on $\cA//\cG$. The difficulty in combining the $7+1$ and $6+1$ dimensional
theories is that we want to find another $3$-form $\rho$ such that we get
the same flow as the {\it gradient flow} of $f_{\rho}$, and this seems to
essentially restrict us to the Calabi-Yau case.   However this restriction
may not be fundamental. We expect that the $(7+1)$ dimensional theory
should be related, from the point of view of the $6$-manifold $Z$, to a  \lq\lq Fukaya category'' of Lagrangian
submanifolds of $\cA//\cG$ which would be, formally, something defined just
by the symplectic structure of $\cA//\cG$.

As usual, there  is a corresponding discussion in the submanifold case. The infinite
dimensional space of symplectic surfaces $\Sigma\subset Z$ has a natural
symplectic form, arising formally as a symplectic quotient \cite{kn:D1}.
Then we get  locally-defined functions on this space by integrating closed
$3$-forms over $3$-dimensional cobordisms.

 \
 
 There is a variant of this discussion which yields deformation of the Special Lagrangian equations in a Calabi-Yau $3$-fold. To explain this we set up some notation. Suppose $L$ is a submanifold of a manifold $M$ and $\psi$ is a $p$-form on $M$. If the restriction of $\psi$ to $L$ vanishes then $\psi$ defines a section of $\Omega^{p-1}(L, N^{*})$ where $N$ is the normal bundle $TM/TL$ of $L$ in $M$. We will denote this section by $\psi_{N}$. Now consider a closed $3$-form $\rho$ and a closed $4$-form $\tau$ on a $6$-manifold $Z$.
There is a locally defined functional $f_{\tau}$ on the space of $3$-dimensional submanifolds of $Z$, defined by integrating $\tau$ over $4$-dimensional cobordisms.
Let ${\cal C}_{\rho}$ be the set of submanifolds $P^{3}\subset Z^{6}$ such that the restriction of $\rho$ to $P$ vanishes. 
We consider the critical points of $f_{\tau}$ {\it restricted to} ${\cal C}_{\rho}$. For any submanifold $P\in {\cal C}_{\rho}$ the restrictions $\rho\vert_{P}, \tau\vert_{P}$ both vanish, the first by definition and the second for dimensional reasons, so we have well-defined bundle-valued forms $\rho_{N}\in \Omega^{2}(N^{*}), \tau_{N}\in \Omega^{3}(N^{*})$. The Euler-Lagrange equation defining the critical points involves a function $f$ on $P$ (which appears as a \lq\lq Lagrange multiplier'') and takes the form
     $$   \tau_{N}= df \wedge \rho_{N}. $$
     Thus we have a system of equations for a pair $(P,f)$:
     $$  \rho \vert_{P}=0 \  , \ \tau_{N}=df\wedge \rho_{N}, $$
     which are analogous to (28). (Of course we need to factor out the constant functions $f$.) When $Z$ is a genuine Calabi-Yau manifold and $\rho,\tau$ are the standard forms the solutions are special Lagrangian submanifolds, as in \cite{kn:RT2}, with $f=0$. More generally if we write $\tau=\omega^{2}$, where $\omega$ is not necessarily closed we can identify these pairs with associative submanifolds in the tube $Z\times \bR$. One would expect that, for generic choices, the solutions $(P,f)$ are isolated and this deformation could be seen as removing a degeneracy in the special Lagrangian equations (which forces the latter to have solutions in moduli spaces of various dimensions, given by the first Betti number of $P$).

\section{Gauge theory on tamed almost-$G_{2}$-manifolds with tubular ends}
\subsection{Compact $7$-manifolds}

We will now make a slightly more detailed analysis of the $6+1$-dimensional
theory. We will do this in the gauge theory setting, but a similar discussion
applies for the submanifold case. We will also restrict attention to a case
when the cross-sections of the ends are genuine Calabi-Yau manifolds.

Suppose that $Y$ is a compact $7$-manifold with an almost $G_{2}$-structure
defined by a $4$-form $\sigma$. 
Let $A$ be any connection on a bundle $E\rightarrow Y$ and form the sequence
of operators
\begin{equation} \Omega^{0}(\frg_{E})\rightarrow \Omega^{1}(\frg_{E})\rightarrow\Omega^{6}(\frg_{E})\rightarrow
\Omega^{7}(\frg_{E}), \end{equation}
as in (21). This is  not in general a complex but we can make a single
operator out of it in the usual way. Use the metric to identify $\Omega^{p}$
with $\Omega^{7-p}$ so we have
\begin{equation}  D_{A}: \Omega^{0}\oplus \Omega^{1} \rightarrow \Omega^{0}\oplus
\Omega^{1}. \end{equation}
 The elliptic operator $D_{A}$ is  self-adjoint if and only if $\sigma$ is closed, but in any case its
symbol is self-adjoint, so the index   is zero.  When $\sigma$ is closed and $A$ is a solution of the $G_{2}$-instanton
equation we get an elliptic complex of Euler characteristic zero and we call
$A$ {\it regular} if the cohomology of this complex vanishes. This implies
in particular that $A$ is isolated in the moduli space of $G_{2}$-instantons.

We will discuss briefly two more technical issues {\it reducible conections} and {\it orientations}. Suppose first that $Y$ is manifold with holonomy equal to $G_{2}$. Then we know \cite{kn:Joyce1} that the harmonic $2$-forms all lie in the
$\Lambda^{2}_{14}$ component in (6). This means that any complex line bundle $L$ over $Y$ admits a $G_{2}$-instanton connection, and in particular such a connection appears as a reducible solutions $A_{0}$ on the bundle $E=L\oplus L^{-1}$. The situation is in some respects similar to that for instantons over a $4$-manifold with negative definite intersection form. The bundle $\frg_{E}$ splits as $\bR\oplus L^{2}$ and, at $A_{0}$, the complex (30) splits into a corresponding sum, with the  interesting part given by 
$$ \Omega^{0}(L^{2}) \rightarrow \Omega^{1}(L^{2})\rightarrow\Omega^{6}(L^{2})\rightarrow
\Omega^{7}(L^{2}). $$
If the cohomology $H^{1}(L^{2}, A_{0})$ vanishes then $A_{0}$ is isolated from irreducible solutions, and this is true also in families of small deformations of the $G_{2}$-structure. Thus, in this case, the irreducible solutions will not affect the enumerative discussion, counting the irreducible solutions. Since the complex has Euler characteristic $0$ we expect that generically, in a family of $G_{2}$-structures,  the cohomology will vanish. However we also expect that it will be non-trivial on some lower dimensional subset. The crucial point however is that the complex structure on $L^{2}$ means that we are considering families of {\it complex linear} operators and for these the generic picture is that cohomology will appear in real codimension $2$. In that case there will be no interaction between the reducible and irreducible solutions in generic $1$-parameter families, which is what are relevant to our purpose. So, granted that a more detailed analysis is necessary we take this as an indication that we can  ignore the potential complications from reducible solutions, and for simplicity we just ignore these reducibles in what follows. (Similar remarks appply to reductions
$S(U(p)\times U(q))\subset SU(l)$.)

Now it seems very reasonable to assume that for   generic $\sigma$ all solutions are regular
Even if this is not the case we could contemplate introducing further, more
artificial,
perturbations of the equations or ideas of \lq\lq virtual cycles'', but let
us assume that perturbations of $\sigma$ suffice.  Our basic goal,  
when the structure is tamed by some closed form $\phi$ is to define a number by counting the solutions  with appropriate signs, and it is the issue of these signs which we take up next. As usual, we mimic the standard discussion in the Casson-Floer theory over $3$-manifolds. We seek to define a \lq\lq relative sign'' $\epsilon(A,A')\in \{\pm 1\}$ for pairs of solutions with the property that 
$$ \epsilon(A,A'')= \epsilon(A,A')\epsilon(A', A''). $$
This gives a way to attach signs to each solution, up to a single overall sign ambiguity. We define $\epsilon(A,A')$ using the spectral flow of  a family of operators
$D_{A_{t}}$, where $A_{t}$ is a path from $A$ to $A'$. Given a path, this spectral flow yields an integer and we set $\epsilon$ to be $1$ or $-1$ as the spectral flow is even or odd. Then the essential thing is to check that this independent of the path, which is the same as saying that the spectral flow around a closed loop (in the space of connections modulo gauge equivalence) is even. Such a loop yields a connection on a bundle $\bE$ over $X=Y\times S^{1}$ and the spectral flow appears as the index of an elliptic operator over $X$. In fact this operator is just the operator apearing in the linearisation of the $\Spin(7)$-instanton equation and can be identified simply as the Dirac operator over $X$, coupled to the bundle $\frg_{\bE}$. Here the spin spaces in $8$ dimensions are regarded as $8$-dimensional real vector bundles. So this question of \lq\lq orientability'' in our $7$-dimensional set-up reduces to an algebro-topological question of showing that the index of this such a coupled Dirac operator over $X$ is {\it even}. (More generally, if the index of all such bundles is divisible by some integer $k$ then we expect the putative Floer theory associated to $Y$ to be $\bZ/k$-graded.)

We can use the Atiyah-Singer index theorem to express this question in terms of characteristic classes. (An explicit formula is given by Lewis in \cite{kn:CL}.)
But we can also give an  argument which avoids detailed calculation. Since it is odd-dimensional, the manifold $Y$ has a nowhere vanishing vector field. This gives a reduction of the structure group of $Y$ to $SU(3)\subset G_{2}$, and hence of $X$ to $SU(3) \subset G_{2} \subset \Spin(7)\subset SO(8)$. With this reduction the spin bundles of $X$ acquire complex structures (corresponding to $V\oplus \bC, V^{*}\oplus \bC$ where $V$ is the fundamental representation of $SU(3)$. Of course we are not saying that these reductions are compatible with the differential geometric structure, but they are at the level of the symbol of the Dirac operator. Thus we can deform the coupled Dirac operator over $X$ to a complex linear operator, and hence the (real) index is even.
However we will not try to develop the theory of orientations any further here.

\

We move on from this brief outline, which indicates how--modulo questions of compactness-- we should define an integer counting the $G_{2}$-instantons on a bundle $E$ over $Y$. 
 For each homology class $b\in H_{3}(Y)$
we consider $SU(l)$ bundles $E$ with $c_{2}(E)$ the Poincar\'e dual of $b$ and (for
simplicity) with $c_{3}(E)=0$. There are at most a finite
number of different topological types and we define an integer $n_{b}$ by
summing the counts above over all such bundles. We expect this to be a deformation
invariant, with respect to perturbations of $\sigma$. 
 The energy bound implies that
$n_{b}=0$ if $[\phi](b)<0$, since the moduli space is then empty. In fact $n_{b}$ vanishes for $b$ outside some proper cone in the half-space $\{b:\phi(b)>0\}$. We package these numbers into a formal series
\begin{equation}   f_{Y}(\psi)= \sum n_{b} \exp(-\langle b,\psi\rangle), \end{equation}
and our further hypothesis is that this converges to yield a holomorphic function of  variable $\psi$ in an open subset of $ H^{3}(Y,\bC)$ containing the points $r [\phi]$, for large enough $r$.

\subsection{$7$-manifolds with tubular ends}

Now we discuss dimension $6$. If we have a Hermitian-Yang-Mills connection on
a holomorphic bundle over a Calabi-Yau manifold $Z^{6}$,
as in the previous section, the linearisation of the equations (28) for pairs $(a,u)$ yields an elliptic deformation complex of Euler characteristic
zero.  We call a solution regular if the cohomology vanishes, and this is
just the same as saying the sheaf cohomology $H^{*}(Z, \frg_{E}\otimes \bC)$
vanishes.  We will assume that all solutions over $Z$ are regular. This is
a very restrictive assumption, and we will return to discuss it further below.

Now suppose that $Y$ is a non-compact manifold with tubular ends each with
a Calabi-Yau cross section, and that we have an
almost $G_{2}$-structure defined by a closed form $\sigma$, compatible with
the product structure on the ends, up to an exponentially decaying term.
 If we fix  a solution over each end then
we have the  notions of an \lq\lq adapted bundle'' and \lq\lq adapted connection''
over $Y$, just as in the usual $3+1$ dimensional theory described in \cite{kn:D2}.
These are, respectively,  a bundle
over $Y$ with a fixed isomorphism with the pull-backs of the chosen bundles
over the ends, and a connection over $Y$ which agrees with the model determined
by the $6$-dimensional solution over each end, up to an exponentially decaying
term. 

Now fix attention on any adapted connection $A$ over $Y$ and form the differential
operator $D_{A}$ as above.  This is formally self-adjoint, just as before.
Over an end, assuming that the connection is actually equal to the model,
we can write $D_{A}$ as
  $$ D_{A}= L \frac{d}{dt}+ Q, $$
  where $L$ is a skew-adjoint algebraic operator  and $Q$ is a first-order
  self-adjoint operator, both over $Z$. Under our hypothesis that $Z$ is Calabi-Yau
  the composite $L^{-1} Q$ is a {\it self-adjoint} first order operator over $Y$ and the 
   analysis of $D_{A}$ follows very closely that in the three-dimensional theory, as described in \cite{kn:D2} for example. We deduce that $D_{A}:L^{2}_{1}\rightarrow
  L^{2} $ is Fredholm
  provided that the operator $Q$ over $Z$ does not have a zero eigenvalue.
  In turn, this is the same as the hypothesis that the solution over $Z$ is
  regular (since the operator $Q$ is essentially the same as that arising
  in the deformation theory over $Z$). 
  
  In the usual Floer theory we go on to consider the Fredholm index, which
  yields an invariant of an adapted bundle. The distinctive feature of the
  theory we are considering here is that the analogous index {\it always  vanishes}. This is just a consequence of the fact, from the general theory,
  that the cokernel of
 $D_{A}$ is represented by the kernel of the formal adjoint, which is the
 same as $D_{A}$. The consequence is that, under our restrictive hypotheses,
  the study of $G_{2}$-instantons
 on a  fixed adapted bundle over $Y$ behaves just like the compact case.
 We define regular solutions in just the same way, we expect that for generic
 $\sigma$ all solutions are regular and that a count of solutions yields
 a deformation invariant (with respect to compactly supported variations
 in $\sigma$)  of the adapted bundle.

 Now we discuss \lq\lq neck stretching sequences'' and gluing constructions.
 For simplicity consider a pair of manifolds $Y_{1}, Y_{2}$ each with one
 end having the same model $Z$  and appearing as the limit of a sequence
 of structures $\sigma_{T}$ on a compact manifold $Y$. (But note that when interchanging $Y_{1},
 Y_{2}$ we have to change the sign of the $3$-form  on  $Z$.) Given
 regular
 $G_{2}$-instantons on adapted bundles over $Y_{1}, Y_{2}$ we wish to construct
 a glued solution over $(Y,\sigma_{T})$ for large $T$. The proof follows
 the familiar Floer theory case closely, with one extra step. For large $T$
 we construct an approximate solution $A_{0}$ with all norms of the error
 term $\sigma_{T}\wedge F(A_{0})$ bounded by decaying exponential functions
 of $T$. Now we seek to solve the equation
     $$  \sigma_{T} \wedge F(A_{0}+\alpha)=  *d_{A_{0}+\alpha}u, $$
     over $Y_{T}$, 
     for a bundle valued $1$-form $\alpha$ and $0$-form $u$. We also impose
     the gauge fixing condition $d^{*}_{A_{0}} \alpha =0$. Schematically, these equations
     can be written as  $$D_{A_{0}} s = s*s + \sigma_{T}\wedge F(A_{0}), $$
     where $s$ is the pair $(\alpha,u)$, the notation $s*s$ denotes a quadratic algebraic
     term and  $D_{A_{0}}$ is our basic elliptic operator. Mimicking the arguments
     in the $3+1$ dimensional case, we get a bound on the operator norm of
     the inverse of $D_{A_{0}}$ which is independent of $T$. Then the inverse
     function theorem shows that when $T$ is large (so $\sigma_{T}\wedge
     F(A_{0})$ is small) there is a small solution $s$. Now the final extra step
     is to observe that, since $\sigma_{T}$ is closed, the Bianchi identity
     implies that in fact $d_{A}u=0$, as we have seen before.
      Hence $A$ is the desired $G_{2}$-instanton.
     
     Still following the familiar pattern, we hypothesise  that when $T$ is large,
     all $G_{2}$-instantons over $(Y,\sigma_{T})$ arise by this gluing construction. Given
     a class $b\in H_{3}(Y)$ we let $c\in H_{2}(Z)$ be the image of $b$
     under the boundary map of the Mayer-Vietoris sequence of $Y=Y_{1}\cup Y_{2}$.
     This is the same as the Poincare dual on $Z$ of the restriction of the
     Poincar\'e dual, $PD_{Y}(b)$, of $b$ on $Y$. We use the symbol $\Sheaf$
     to denote a solution over $Z$ on a bundle with Chern class the Poincar\'e
     dual of $c$.  Then our hypothesis gives a gluing formula of the shape
     \begin{equation}  n_{b}= \sum_{\Sheaf} \sum_{ E_{1}, E_{2}} n(E_{1})n(E_{2}), \end{equation}
     where  in the inner sum $E_{1}, E_{2}$ run over adapted bundles with common limit $\Sheaf$
     and such that the glued bundle over $Y$ has $c_{2}=P.D. (b)$. Of course
     we have similar formulae when we glue manifolds with more than one end.
     
     It is important to emphasise that if we work with almost $G_{2}$ structures
     where $\sigma$ is not closed we would get a different, richer theory, more
     like Floer theory over $3$-manifolds and $4$-manifolds with tubular
     ends.  The operator $D_{A}$ is not self-adjoint
     and the index gives a non-trivial invariant of adapted bundles. Then
     we would get moduli spaces of different dimensions, depending on the
     index. However we will not pursue this further here.

     \subsection{Holomorphic bundles over moduli space}
     
     We have now reached the main question we wish to address in this paper.
     The gluing formula (33) is, in the general Floer theory framework, at
     the \lq\lq chain level''. As we vary the Calabi-Yau structure  on $Z$ the
     solutions $\Sheaf$ vary and we do not have a canonical way to identify
     them. Further, even if choose such an identification locally the individual
     numbers $n_{E_{1}}, n_{E_{2}}$ will change. So we seek a more invariant
     way of expressing the formula much as, in the ordinary Floer theory,
     we pass from the \lq\lq chains'' to the Floer homology groups.

   Let us for simplicity suppose that the inclusion of $Z$ in $Y$ induces
   an injection on $H_{3}$ and fix a coset $[b_{0}]$ in $H_{3}(Y)/H_{3}(Z)$. We restrict attention to classes $b$ in this coset, which have the same image
$c$ in $H_{2}(Z)$. Now we consider a series like (32)
\begin{equation}  g_{Y}(\psi)= \sum_{b\in H_{3}(Z)} n_{b+b_{0}} \exp(-\langle b , \psi\rangle), \end{equation}
where now $\psi$ lies in an suitable open set in $H^{3}(Z,\bC)$, containing  the points $r[\theta]$ for large enough $r$.

Fix the class $[\omega]\in H^{2}(Z)$ and consider the moduli space $\cM$ of pairs $(I,\theta)$ where $I$ is a complex structure on $Z$ which
admits Kahler metrics in the class $[\omega]$ and $\theta$ is a nowhere-zero holomorphic $3$-form. By the Torelli theorem for Calabi-Yau manifolds this is a quotient $\tcM/\Gamma$ where $\Gamma$ is the symplectic mapping class group and
$\tcM$ is immersed in $H^{3}(Z,\bC)$. We have an obvious $\bC^{*}$-action on $\cM$, multiplying the holomorphic form by a constant. We can also define a \lq\lq norm'' function by
$$  \Vert (I, \theta) \Vert^{2}= \int_{Z} \theta\wedge \overline{\theta}. $$
and for $R>0$ we write $\cM_{R}$ for the points in $\cM$ of norm greater than $R$.  We suppose $R$ is some fixed, suitably large, number. Then the restriction of the holomorphic function $g_{Y}$ defines a holomorphic function on an open set in $\cM_{R}$ containing the ray
$l_{\theta}=\{ I,  r\theta): r\geq 1\}$. We denote this restriction also by $g_{Y}$.

Now we make the assumption that {\it for generic points in the moduli space
$\cM$ all solutions $\Sheaf$ are regular}.

\begin{formal}
Under this assumption:
\begin{itemize}
\item  There should be a holomorphic vector bundle $\cE\rightarrow \cM_{R}$, associated to the class $c\in H_{2}(Z)$, and a canonical isomorphism $(-1)^{*} \cE \cong \cE^{*}$. \item There should be an invariant $g_{Y_{1}}$ which is a holomorphic section of $\cE$ over a neighbourhood of the ray $l_{\theta}$. Likewise
there should be  an invariant $g_{Y_{2}}$ which is a holomorphic section of $(-1)^{*}\cE$ over a neighbourhood of $l_{\theta}$.
\item The function $g_{Y}$ should be the dual pairing $\langle g_{Y_{1}}, g_{Y_{2}}\rangle$, formed using the isomorphism above \end{itemize}
\end{formal}
In other words, in this situation, the structure analogous to the Floer homology of a 3-manifold is the holomorphic vector bundle $\cE$ over the moduli space, and its sheaf
of holomorphic sections.

\

{\bf Remarks}

\

1. The mapping class group may not act freely on $\tcM$,  in which
case $\cM$ is an orbifold. Then we more should strictly work with an \lq\lq
orbi-bundle''
over $\cM$. 

2. When we restrict the function $g_{Y}$ from an open set in $H^{3}(Z,\bC)$
to
$\tcM$ we could lose information. But our bundles actually extend over a
thickening of $\cM$ obtained from a quotient of such an open set and if formulated
this way we get a gluing formula for the original function $g_{Y}$.

 \
 
Now we will explain the construction of the bundle.
Suppose  that $\theta,\theta'$ are two nearby  Calabi-Yau structures on $Z$
and that they are generic in the sense above, so all solutions  are regular.  We have a collection of critical points $\{\Sheaf\}$ associated to  $\theta$ (with the fixed class $c$) and another collection $\Sheaf'$ associated to $\theta'$.  Let $V,V'$ be the complex vector space with basis elements $\langle\Sheaf\rangle$ and $\langle\Sheaf'\rangle$ respectively.
Write $\langle\Sheaf\rangle^{*}$ for elements of the dual basis of $V^{*}$.  We can choose  a tamed almost $G_{2}$-structure on the topological cylinder $\bR\times Z$ which is asymptotic to that defined by $\theta$ at $-\infty$ and by $\theta'$ at $+\infty$. Fix, for the moment, closed $3$-forms $\upsi$ representing cohomology classes $\psi\in H^{3}(Z)$ and extend these to closed $3$-forms
$\opsi$ on the topological cylinder, compatible with the product structure on the ends. We choose the taming form $\phi$ as the  representative of ${\bf Re}(\theta)$. Then for each adapted bundle $E$ over the cylinder we have the
following

\begin{enumerate}
\item Asymptotic  limits $\Sheaf(E) \in \{ \Sheaf\},\Sheaf'(E)\in \{\Sheaf'\}$, at $t=\pm
\infty$ respectively.
\item  A real number
$$   \upsi(E)= -\frac{1}{8\pi^{2}}\int_{Z\times\bR} \opsi \ \Tr(F_{A})^{2},
$$ where $A$ is any adapted connection. This is independent of the choice
of connection and of the representative $\opsi$, given a choice of $\upsi$.
\item An integer $n(E)$ counting the number of $G_{2}$-instantons. 
\end{enumerate}

Now we define a holomorphic function with values in ${\rm Hom}(V,V')$
\begin{equation}   G_{\theta,\theta'}(\psi)= \sum_{E} n(E) \exp(-\upsi(E)) \langle\Sheaf(E)\rangle^{*}\otimes \langle\Sheaf'(E)\rangle.\end{equation}
This will be defined in a neighbourhood of the ray $l_{\theta} $. If we change the
choice of representive $\upsi$ of a class $\psi$ to $\upsi+ d\lambda$ we
change the numbers $\upsi(E)$ to
$$ \upsi(E)+ \int_{Z} \left(\Tr(F_{\Sheaf(E)}^{2})- \Tr(F_{\Sheaf'(E)}^{2}\right) \lambda. $$
This has the effect of changing the linear map $G_{\theta, \theta'}(\psi)$ to 
\begin{equation}\Lambda_{\theta}(\psi)
G_{\theta, \theta'}(\psi)\Lambda_{\theta'}(\psi)^{-1}, \end{equation}
where in our bases $\Lambda_{\theta}, \Lambda_{\theta'}$ are diagonal matrices
with entries given by the exponentials of 
$$  \int_{Z} \Tr F_{\Sheaf}^{2}\ \lambda \ \ , \int_{Z} \Tr F_{\Sheaf'}^{2}
\ \lambda, $$
respectively.
If $\theta''$ is another nearby generic structure and we fix the same representatives
the gluing formula, extended to this situation (when we glue together two
topological cylinders) yields

\begin{equation}  G_{\theta,\theta''}(\psi)= G_{\theta, \theta'}(\psi) \circ G_{\theta', \theta''}(\psi). \end{equation}
Thus to define our holomorphic bundle we decree that for each generic structure $\theta$ and choice of representatives $\upsi$ the bundle has a canonical trivialisation 
over some neighbourhood
$U_{\theta}$ of $l_{\theta}$ in the moduli space. Changing the representatives
$\upsi$ changes the trivialisation by multiplying by the diagonal metric
$\Lambda_{\theta}$. On the
  the overlaps
$U_{\theta}\cap U_{\theta'}$ we use  the maps $G_{\theta,\theta'}(\psi)$ as transition functions. The gluing formula (37) is the cocycle condition giving the consistency of the set of transition functions, and
the isomorphism $(-1)^{*}(\cE)=\cE^{*}$ is induced by the fact that the solutions
$\{\Sheaf\}$
of the equations  defined by the $3$-forms
$\rho, -\rho$ have an obvious identification.

Now to define the section $g_{Y_{1}}$ corresponding to $Y_{1}$ we work initially
in a canonical trivialisation around $\theta$ and with a choice of the representatives
$\upsi$. By our hypothesis these can be extended to closed forms over $Y_{1}$,
compatible with the product structure on the end. Then for an adapted bundle
$E_{1}$ over $Y_{1}$ and a class $\psi$ in $H_{3}(Z)$ we define $\upsi(E_{1})$ in
the same manner as before. We also have a limit $\Sheaf(E_{1})$ and a number $n(E_{1})$
The we set 
$$  \tilde{g}_{1}(\psi)= \sum_{E_{1}} n(E_{1}) \exp(- \upsi(E_{1})) <\Sheaf(E_{1})>, $$
a $V$-valued holomorphic function on a neighbourhood of $l_{\theta}$. The
gluing formula shows first that this yields a well-defined section $g_{1}$
of
 $\cE$---independent of the choice of $\theta$ and of $\upsi$--- and second that $g_{Y}=<g_{Y_{1}}, g_{Y_{2}}>$.

\

{\bf Remarks}

1. Our  bundle should have the property that its rank is equal to the DT invariant. We can think of the rank as the $0$-degree component of the Chern character, and it might be that there is an extension of these ideas to a formula for all of ${\rm ch}\ \cE$. Note that  Thomas' theory discusses a fixed Calabi-Yau manifold and there should be extensions of this which yield cohomology classes in the moduli space $\cM$, using the universal family. 

2. In this Section we have fixed attention on the $(6+1)$-dimensional theory
but it is natural to wonder if there is some yet higher structure, associated
to the $8$-dimensional geometry. Roughly speaking, we would expect this to assign
Floer groups to a compact $G_{2}$-manifold $Y$ and one would like some machinery
to compute these when $Y=Y_{1}\cup_{Z}Y_{2}$. We make one observation in
this direction. Suppose we have a vector bundle over a space $B$ given by transition
functions $g_{\alpha, \beta}$ with respect to a cover $U_{\alpha}$ of $B$.
Thus we have matrix entries $g_{\alpha, \beta, i,j}(z)$ which are functions
of  $z\in U_{\alpha}\cap
U_{\beta}$ and
$$   g_{\alpha, \gamma, i, k}(z)= \sum_{j} g_{\alpha, \beta, i,j} \ g_{\beta,
\gamma, j,k}, $$
on $U_{\alpha}\cap U_{\beta}\cap U_{\gamma}$. Now suppose we have a chain
complex $C^{*}$ and a chain automorphism $T:C^{*}\rightarrow C^{*}$ so we
have a Lefschetz number $L(C^{*},T)=\sum (-1)^{i} \Tr T_{i}$. If we have a family
of such pairs, parametrised by $z$, the Lefschetz number becomes a function
of $z$. So we may envisage a structure given by chain complexes
$C^{*}_{\alpha, \beta, i,j}$, with automorphisms, parametrised by $z\in U_{\alpha}\cap
U_{\beta}$, such that on the triple overlaps
$$   C_{\alpha,\gamma,i,j,k} \sim \bigoplus_{j}  {\rm Hom}(C^{*}_{\alpha,\beta,i,j},
C^{*}_{\beta, \gamma,j,k}), $$ 
where $\sim$ is a suitable equivalence relation, at least as strong as chain
homotopy
equivalence, compatible with the automorphisms. Then the Lefschetz numbers of the chain complexes give the transition
functions of a holomorphic bundle. Possibly there is some structure of this kind on the moduli space
$\cM$ which would refine the bundle $\cE$, in the same way as the Floer homology of $G_{2}$-manifolds should refine  the Casson invariant.

\section{Finite-dimensional analogue} 
\subsection{Morse-Novikov Theory}
We begin with  some basics. Take a manifold $A$, and a class $[\alpha] \in H^1(A, \bC)$.
This gives a local system, which we may think of either as a rank 1
complex vector bundle on $A$ with flat connection $\alpha$, or as a representation
(the monodromy of the connection)
$$\rho:\pi_1(A)\to \bC^*$$
$$\rho(\gamma) =  \exp ( \int_\gamma \alpha)$$
How do the cohomology groups $H^*(A, \rho)$ of this local system behave as we vary $[\alpha]$? The
answer is that they form the fibres of a holomorphic sheaf $\mathcal{E}$
over the vector
space $H:= H^1(A, \bC)$. We show this as follows. An element $\gamma\in \pi_1(A)$
gives a linear function on $H$. So we can let $\pi_1(A)$ act on the
ring of holomorphic functions on $H$ by
$$ \gamma: f \mapsto e^\gamma f$$
This is an action by module automorphisms, so we have a local system on $A$ whose fibre is the rank 1 free $\mathcal{O}_H$-module. The cohomology of
this local system $\mathcal{E}:=H^*(A, \mathcal{O}_H)$ is then also an $\mathcal{O}_H$-module,
i.e. a sheaf on $H$, and its fibre at the point $[\alpha]$ is the cohomology
of the local system associated to $[\alpha]$.

We will now review Morse-Novikov theory in finite-dimensions, more specifically the reformulation
due to Burghelea and Haller \cite{kn:BH}. Let $\Theta$ be a closed real 1-form
on a finite-dimensional Riemannian manifold $A$, and assume that $\Theta$ is \lq\lq Morse'', i.e. locally
the differential of a Morse function. Then the pull-back of $\Theta$ to the
universal cover $\tilde{A}$ is globally the differential of a Morse function,
and we may form the cell-complex $C_\bullet$ given by the unstable manifolds of the
zeroes of $\Theta$ on $\tilde{A}$ in the usual way. This cell-complex is
obviously periodic with
respect the action of $\pi_1(A)$.

Now suppose we have a rank-1 complex local system on $A$ given by the class
of some closed complex 1-form $\alpha$. Then we can use our periodic cell-complex $C_\bullet$ in $\tilde{A}$ to produce a chain-complex that in good cases will calculate the cohomology of this local
system. What we do is twist the differential on $C_\bullet$ by $\alpha$ and then take $\pi_1(A)$ invariants. Assuming that $\Theta$ is in fact Morse-Smale for the metric on $\tilde{A}$, we can describe this chain-complex explicitly in the following way: it  has a basis given by the set $Z_A(\theta)=\{p\}$ of zeroes of $\Theta$ in $A$, and differential
\begin{equation}\partial_{\alpha} p = \sum_{i(q)=i(p)-1}\; \sum_{\gamma:p\to q} e^{\int_\gamma \alpha}q \end{equation}
where $i(p)$ denotes the index of $p$, and $\gamma$ runs over all flow-lines of $\theta$ connecting $p$ and $q$. The differential clearly varies holomorphically
in $\alpha$, so the homology of the chain complex is a holomorphic sheaf
$\mathcal{E}$
on the space $\Omega_c^1$
of closed complex 1-forms. If we choose a section $H \to \Omega_c^1$ (e.g.
by using Hodge theory) we can pull-back and get a holomorphic sheaf
on $H$. The fibre of this sheaf at $[\alpha]$ is the cohomology of the local system on $A$ given by $[\alpha]$.

There are two problems with this description. Firstly, since $\tilde{A}$ is non-compact the cell-complex $C_\bullet$ may not cover the whole of $\tilde{A}$ (indeed $\theta$ may have no zeroes at all), so we can't guarantee that we are genuinely calculating the cohomology of the local system $[\alpha]$. Secondly, the sum (37) may not converge. It is shown in \cite{kn:BH}, that neither of these issues is a problem when $\alpha = s\theta$ for $s$ a complex number with $\Re(s)>\rho$, where $\rho\in [0,\infty]$ is some numerical invariant depending on $\theta$ and the metric.
For this range of $\alpha$, the components of $\partial_\alpha$ do converge,
and the chain-complex does calculate the cohomology of $[\alpha]$ correctly.
The invariant $\rho$ is conjectured to be always finite, and in many examples
is zero.

Let's suppose something stronger: that we get convergence and the correct cohomology in an open neighbourhood of the set $\Re(\alpha)=\theta$. This includes the set $\{\alpha=(1+i\lambda)\theta,
\lambda\in \bR\}$, where we do indeed get convergence and the right cohomology
at least when $\rho<1$. If this holds, for any generic $\alpha_0 \in \Omega^1_c$
we have an analytic way of constructing our sheaf $\mathcal{E}$ in a neighbourhood
of $\alpha_0$. We set $\Theta = \Re(\alpha_0)$, then the Morse-Novikov
chain-complex (38) constructed from $\theta$ will be valid in a neighbourhood of $\alpha_0$, so we may define $\mathcal{E}$ to be its homology.

Now consider the special situation in which $A$ is a complex manifold and $\Theta$ is the real part of a holomorphic $1$-form. Then the indices of all critical points are equal and, in a generic situation there are no flow lines between critical points. For any fixed $\alpha$ the cohomology $H^{*}(A,[\alpha])$ is simply a vector space with one basis element for each zero of $\Theta$.
The sheaf on $\mathcal{E}$ on $H^{1}(A,\bC)$ is locally free and we write
$\mathcal{E}$ also for the corresponding holomorphic vector bundle.  

We wish to apply this construction, formally, to the infinite-dimensional manifolds $\mathcal{X}$
arising from a Calabi-Yau
threefold $Z$. Thus $\mathcal{X}$ is either the space of unitary connections
with $F\wedge \omega^{2}=0$, modulo gauge equivalence, or the space of symplectic surfaces $\Sigma \subset Z$. Assume for simplicity that we are working with genuine Calabi-Yau structures. Then given a point  $(I,  \theta)\in \tcM$ we get an induced complex structure on $\mathcal{X}$ making it a Kahler manifold, and the Chern-Simons construction gives a holomorphic $1$-form on $\mathcal{X}$ whose zeros correspond to stable holomorphic bundles or complex curves respectively.
More generally, we have a map from $H_{3}(Z)$ to $H^{1}(\mathcal{X})$ and the Chern-Simons construction gives a definite representative of this---a $1$-form on $\mathcal{X}$ corresponding to a closed $3$-form on $Z$.  

Allowing the complex structure to  vary, we get a map from $\tilde{\mathcal{M}} \to H^{1}(\cX;\bC)$, which
is in fact a holomorphic embedding. So we should get a sheaf over $\tilde{\mathcal{M}}$ as the  restriction) of the twisted cohomology sheaf $\mathcal{E}$, assuming of course that we had made sense of the latter. Then everything is invariant under the group $\Gamma$ so we can descend to the moduli space $\mathcal{M}$.

It should now be fairly clear how our (conjectural) construction fits into this picture. Although we do not have any reason to believe that there is a cohomology sheaf over the whole of $H^{3}(Z,\bC)$ we can make sense of this over a neighbourhood of $\tcM$. An individual fibre, over a generic point, is rather uninteresting--having an almost-canonical basis. The interesting structure appears in the way these are fitted together into a vector bundle.
To understand this we go back to our  finite-dimensional situation and
consider Morse theory on the universal cover.

Let $\theta$ again be a closed real 1-form on a finite-dimensional Riemannian manifold $A$,
and (by abuse of notation) let $\theta$ also denote its pull-back to the universal
cover $\tilde{A}$. Assume that $\theta$ is Morse-Smale, then we want to consider the resulting Morse complex for $\tilde{A}$.

Choose lifts $\tilde{p}$ of each zero $p\in A$ of $\theta$,
so the zeroes of $\theta$ on $\tilde{A}$ are 
$$Z_{\tilde{A}}(\theta)=\{\gamma(\tilde{p}),
\gamma \in \pi_1(A)\}$$
This is a basis for the Morse complex. The set of flow-lines is $\pi_1(A)$-invariant, so when $i(q)=i(p)-1$
we may define $N^{\gamma^{-1}\phi}_{pq}$ to be the (signed) number of flow-lines between $\gamma(\tilde{p})$
and $\phi(\tilde{q})$, this depends only the product $\gamma^{-1}\phi$. Then the differential in the Morse complex is
$$ d(\gamma(\tilde{p})) = \sum_{i(q)=i(p)-1}\;\sum_{\phi\in \pi_1(A)} N^{\gamma^{-1}\phi}_{pq}
\phi(\tilde{q})$$ 

Since Morse homology is actually isomorphic to singular homology, it doesn't
change as we vary  $\theta$ and the metric.  Floer  observed in $\cite{kn:F}$ that
we can give an \textit{a priori} proof of this fact, without reference to
singular homology, by counting flow-lines in
families. Suppose we have a suitably generic path $(\theta_t, g_t)$ of 
($\pi_1(A)$-invariant)
closed 1-forms and metrics, where $(\theta_0, g_0)$ and $(\theta_1, g_1)$
are Morse-Smale. Then using some standard function
we can define a 1-form on $\tilde{A}\times[0,1]$ which has zeroes 
$$\{\theta_0=0\}\times\{0\} \bigsqcup \{\theta_1=0\}\times\{1\} $$
Let $M^{\gamma^{-1}\phi}_{pq}$ be the (signed) number of flow-lines 
$$\gamma(\tilde{p})\times \{0\} \longrightarrow \phi(\tilde{q})\times\{1\}$$
 where $i(p)=i(q)$, and define a linear map $\Psi$
between the Morse complexes for $(\theta_0, g_0)$ and $(\theta_1, g_1)$ by
\begin{equation}\label{Psi} \Psi(\gamma(\tilde{p})) = \sum_{i(q)=i(p)}\; \sum_{\phi\in \pi_1(A)} M^{\gamma^{-1}\phi}_{pq}
\phi(\tilde{q})\end{equation}
Then it can be shown that $\Psi$ is a map of chain complexes, and is functorial
with respect to composition of paths. 
Furthermore, if we have a homotopy
between two such paths, then by a similar construction there is an induced homotopy between the corresponding
chain maps. It follows that $\Psi$ is a homotopy equivalence. (Here we are assuming, in the finite-dimensional case, that we do not run into problems due to the noncompactness of $\tilde{A}$.)

Now suppose that $\theta_0=\Re(\alpha_0)$ and $\theta_1=\Re(\alpha_1)$ for some complex 1-forms $\alpha_0$ and $\alpha_1$. We can use the above discussion to understand how the Morse-Novikov
complexes for $(\alpha_0, g_0)$ and $(\alpha_1, g_1)$ are related. To do this, we need the
following construction.

Let $(\theta_0, g_0)$ and $(\theta_1, g_1)$ be any two Morse-Smale pairs.
Suppose, we have a $\pi_1(A)$-invariant linear map 
$$T: \left<Z_{\tilde{A}}(\theta_0)\right>_\bR \to \left<Z_{\tilde{A}}(\theta_1)\right>_\bR $$ of the vector spaces underlying the corresponding Morse complexes. Then
$$T: \gamma(\tilde{p}) \mapsto \sum_{q, \phi} T^{\phi^{-1}\gamma}_{pq}
\phi(\tilde{q})$$
for some set of  real numbers $\{T^\gamma_{pq}\}$. Examples are the Morse
differential $d$ when $(\theta_0, g_0)=(\theta_1, g_1)$, and the map $\Psi$ (\ref{Psi}) when $(\theta_0, g_0)\neq (\theta_1, g_1)$.  Then for any closed
complex 1-form $\alpha$ on $A$ we
may formally define a map
$$\hat{T}: \left<Z_{A}(\theta_0)\right>_\bC \to \left<Z_{A}(\theta_1)\right>_\bC
$$
by setting
$$ \hat{T}_\alpha(p) = \sum_{q,\gamma} T^\gamma_{pq} e^{\int_{\tilde{p}}^{\gamma(\tilde{q})}
\alpha} q $$
though of course this may not converge. It is easy to check that this
process is (formally) a homomorphism, i.e. if we have three Morse-Smale pairs
$(\theta_0, g_0),(\theta_1, g_1),(\theta_2, g_2)$ and maps 
$$\left<Z_{\tilde{A}}(\theta_0)\right>_\bR \stackrel{T}{\to} \left<Z_{\tilde{A}}(\theta_1)\right>_\bR
\stackrel{U}{\to} \left<Z_{\tilde{A}}(\theta_2)\right>_\bR $$
then 
\begin{equation}\label{hathomo}\hat{T}_\alpha\hat{U}_\alpha=\widehat{TU}_\alpha
\end{equation}
Also applying this to the Morse differential we get 
$$ \hat{d}_\alpha(p) = \sum_{i(q)=i(p)-1}\;\sum_{\gamma} N^\gamma_{pq} e^{\int_{\tilde{p}}^{\gamma(\tilde{q})} \alpha}
 q $$
which is just the Morse-Novikov differential $\partial_\alpha$.

Now take two generic pairs $(\alpha_0, g_0)$, $(\alpha_1, g_1)$ of closed complex 1-forms and metrics on $A$, and choose a generic path btween
them in the space of closed complex 1-forms and metrics. Take the real parts of all the 1-forms to get the same data in the space of closed real 1-forms and metrics. We get two Morse complexes with differentials $d^0$ and $d^1$,
and a chain map $\Psi$ between them. Assume there is some neighbourhood $\mathcal{U} \subset \Omega^1_c$, containing
$\alpha_0$ and $\alpha_1$, such that the maps $\hat{d}^0_\alpha, \hat{d}^1_\alpha$
and $\hat{\Psi}_\alpha$ converge for all $\alpha\in \mathcal{U}$. Then over $\mathcal{U}$ we have two  Morse-Novikov
complexes of holomorphic vector bundles on $\mathcal{U}$, given by the differentials
$\hat{d}^0$ and $\hat{d}^1$.  However, by (\ref{hathomo}) we know that $\hat{\Psi}$
is a holomorphic chain map between the two complexes.

If we have two such paths and a homotopy between them, we get two chain maps
$\Psi^1$ and $\Psi^2$ between the Morse complexes, and chain-homotopy $\Xi$.
Then, assuming everything converges in $\mathcal{U}$, we get two holomorphic
chain maps $\hat{\Psi}^1$ and $\hat{\Psi}^2$ between our Morse-Novikov complexes
on $\mathcal{U}$ and a holomorphic homotopy $\hat{\Xi}$ between them. Similarly,
assuming that all the necessary maps converge in $\mathcal{U}$, we can show
that our two Morse-Novikov complexes are homotopy equivalent. 

In the model situation we are considering, the complex 1-forms on $A$ that
are actually holomorphic, so at generic points the Morse-Novikov chain complexes
have no differential.  This means the homotopy equivalences $\hat{\Psi}$ are just isomorphisms, and that if $\hat{\Psi}^1$ and $\hat{\Psi}^2$ are homotopic maps then they
are actually identical. 

We now have  a picture of how our bundle $\mathcal{E}$ on $H^{1}(A,\bC)$
is built up.  Suppose we have a simply-connected region $\mathcal{U}\subset\mathcal{M}$
consisting of only two \lq\lq chambers'' $\mathcal{U}_1$ and $\mathcal{U}_2$ in which
the Morse-Smale condition holds, separated by a wall on which it fails. 
Firstly, take two points $x,y \in \mathcal{U}_1$, and a path between them.
The corresponding map $\Psi$ counts flow-lines in the family of 1-forms and
metrics given by the path, but since the Morse-Smale condition holds everywhere
along the path we will just see a single flow-line from any zero of $\alpha_x$
to the corresponding zero of $\alpha_y$, so $\Psi$ is just the identity map.
This means that $\mathcal{E}$ should be trivial over each region $\mathcal{U}_1$
and $\mathcal{U}_2$.
However, if we have a path between $x\in \mathcal{U}_1$ and $y\in \mathcal{U}_2$ then the
family will contain additional flow-lines arising from the \lq\lq unexpected''
flow-lines that occur on the wall, so the map $\Psi$ will be non-trivial.
Passing to the holomorphic version, we get the transition function $\hat{\Psi}$
that we should use to patch $\mathcal{E}$ when we cross the wall.

We can formulate this construction a little more precisely. Pick any $x\in
\mathcal{U}_1$ and $y\in \mathcal{U}_2$ and any path between them. The corresponding
map $\Psi_{xy}$ is independent of the choice of path since $\mathcal{U}$
is simply-connected. Assume that $\hat{\Psi}_{xy}$
converges in all of $\mathcal{U}$. Let $V_x, V_y$ be the trivial bundles on $\mathcal{U}$
with fibres $\left<Z_{\mathcal{A}}(\alpha_x)\right>_\bC$ and 
$\left<Z_{\mathcal{A}}(\alpha_y)\right>_\bC$ respectively. Take another copy
of each of these bundles, and denote them $V_{xy}:=V_x$ and $V_{yx}:=V_y$.
Define
$$\mathcal{E}_\mathcal{U} = {\rm coker} ( V_{xy}\oplus V_{yx} \stackrel{\Phi }{
\longrightarrow} V_x\oplus V_y)$$
where
$$ \Phi = \left( \begin{array}{cc} {\rm id}& \hat{\Psi}_{xy} \\
\hat{\Psi}^{-1}_{xy} & {\rm id }\end{array} 
\right)$$

Then over $\mathcal{U}_1$ we have $\mathcal{E}_\mathcal{U} \cong V_x$ by projecting onto the first
factor, and over $\mathcal{U}_2$ we have $\mathcal{E}_\mathcal{U} \cong V_y$ by projecting onto the second factor, but across the wall these isomorphisms differ by the transition
function $\hat{\Psi}_{xy}$. This is precisely the vector bundle
that we wanted.

In fact this picture of a finite wall-and-chamber decomposition is misleading.
Since there are infinitely many homotopy classes in which non-Morse-Smale
flow-lines can appear, we actually expect the set of walls to be dense in
the parameter space $H^{1}(A,\bC)$. However, this construction adapts easily. Let $\mathcal{U}$
again be simply connected, but have possibly infinitely many walls in it.
For every point $z\in \mathcal{U}$ that doesn't lie on a wall we let 
$V_z$ be the trivial bundle on $\mathcal{U}$ spanned by $Z_{\mathcal{A}}(\alpha_z)$,
and for every distinct ordered pair $x,y$ of such points we take a copy $V_{xy}:=V_x$.
For every pair $x,y$ we get a map $\Psi_{xy}$ of Morse complexes, and again
this is independent of the choice of path between $x$ and $y$. Assume that
$\hat{\Psi}_{xy}$ converges in all of $\mathcal{U}$ for all $x,y\in \mathcal{U}$.
Define
$$ \mathcal{E}_\mathcal{U} = {\rm coker}( \bigoplus_{x, y} V_{xy} \stackrel{\Phi}{\longrightarrow} \bigoplus_z V_z)$$
where the map $\Phi$ has components
$$ \delta_{xz} {\rm id }+ \delta_{yz}\hat{\Psi}_{xy} : V_{xy} \to V_z$$
At any generic point $z$ we have $\mathcal{E}_\mathcal{U} \cong V_z$ by projecting
onto that factor, so in particular $\mathcal{E}_\mathcal{U}$
is a finite rank vector bundle.

\subsection{Discussion}
In our infinite-dimensional situation we interpret Floer's time dependent vector fields, used in constructing the chain homotopy, as almost-$G_{2}$ structures on topological cylinders. The count of flow lines for the time dependent vector field gives the count of $G_{2}$-instantons, and the Novikov series handle the passage to the universal cover. Then it should be clear how our definition of the holomorphic bundle over $\cM$ matches up with the finite-dimensional discussion. To sum up, {\it we interpret the holomorphic bundle $\mathcal{E}$ over $\cM$ as the \lq\lq middle-dimensional'' cohomology of the infinite dimensional space $\cX$, with coefficients in the local system over $\cX$ defined by classes in $H^{3}(Z,\bC)$}.

When Floer first introduced his theory he gave a different treatment of the parameter-dependence of his homology which was more explicit but technically more complicated. Later, he introduced the \lq\lq time-dependent vector field'' trick which gave a simpler approach, but the original method had the advantage of giving a more precise description of the chain homotopies. In a similar way we can give a more explicit description of the co-efficients which define the transition function of our holomorphic bundle. This is closer to the discussion in \cite{kn:DT} and in \cite{kn:DSal1}. The \lq\lq Novikov'' aspect of the theory is confusing here, so let us begin by considering a finite dimensional complex manifold $\tA$ and a proper holomorphic map $f:\tA\rightarrow \bC$ with a finite number of critical values $z_{1},\dots, z_{N}\in \bC$, where we suppose each $z_{i}$   corresponds to a unique critical point $A_{i}\in \tA$. The gradient curves of the function $Re(f)$ on $\tA$ map under $f$ to line segments with fixed imaginary part. So if, as will generically be the case, the imaginary parts of the $z_{i}$ are all different there can be no gradient curves joining the critical points. However we are free to rotate the picture,  multiplying $f$ by any complex number of modulus $1$. Given a pair $z_{i}, z_{j}$ we can choose this phase so that we are allowed to have  gradient curves joining $A_{i}, A_{j}$, mapping to the line segment $\overline{z_{i} z_{j}}$ in $\bC$, which generically will not contain any other $z_{k}$. The count of such gradient curves gives a number $N_{ij}$. This is the intersection number of the vanishing cycles of $f$ at the two points, when the fibre $f^{-1}(t)$ is transported along the segment  $\overline{z_{i}z_{j}}$ so that the two vanishing cycles can be viewed as homology classes in the same space.  Now suppose we have a generic $1$-parameter family of such situations. Then in the family a third point $z_{k}$ may move across the segment $\overline{z_{i}z_{j}}$. When this happens the number $n_{ij}$ changes by $\pm n_{ik} n_{kj}$. Another way of expressing this is that we change the homotopy class of the path in $\bC\setminus\{z_{1},\dots, z_{N}\}$ used to identify the fibres near $z_{i}$ and $z_{j}$. The parallel transport along the two paths differs by monodromy around $z_{k}$: a Dehn twist in the vanishing cycle associated to $z_{k}$. Now given a generic $f$ we let $V_{f}$ by the vector space with basis symbols
$\langle A_{i}\rangle$ associated to the critical points. As we move in a generic $1$-parameter family from $f_{0}$ to $f_{1}$ say the $A_{i}$ move continuously so we have a naive, and rather trivial, isomorphism between $V_{f_{0}}$ and $V_{f_{1}}$ (depending on the path). We modify this in the following way. If in the family the critical point $z_{j}$ moves across the ray $\{z: Im(z)=Im(z_{i}), Re(z)>Re(z_{i})\}$ then we map the basis element $\langle A_{i}\rangle$, before the crossing, to $\langle A_{i}\rangle \pm n_{ij} \langle A_{j}\rangle $, after the crossing.

Now replace the hypothesis that $f$ be proper and with a finite number  of critical points by the situation we had before, where $\tA$ is the universal cover of some $A$, so there are a finite number of critical points up to the action of
$H_{1}(A)$. The same procedure gives a more explicit recipe for the local trivialisation of the twisted cohomology bundle, except now we need to keep track of the relative homotopy class of paths between critical points, which brings in the Novikov series. Translated into our infinite dimensional picture this gives a way to describe the transition functions of our holomorphic bundle, using $G_{2}$-instantons on genuine tubes $\bR\times Z$, exploiting the freedom to multiply the $3$-form $\theta$ by a phase.

\

\

With all of this discussion in place, we return to discuss our main hypothesis,
that for generic points in the moduli space $\cM$ all solutions are regular.
This seems to be a serious restriction, not often satisfied in practice.
However there is an obvious strategy for removing  this restriction. We work
with suitable generic perturbations of the Calabi-Yau structure, involving
triples $\omega, \rho,\rho'$ with $\omega\wedge \rho=0$. It is very reasonable to expect that for generic
perturbations of this kind all solutions are regular.  But as we explained
in Section 3 we have then to give up the assumption that $\sigma$ is closed,
so we get nonzero Fredholm indices for adapted bundles. But this just means
that, in the finite-dimensional analogue, we need to compute twisted cohomology
using a $1$-form with zeros of different indices so we can have a nontrivial
chain complex. What one would expect is that it is possible to define a collection
of sheaves over the moduli space $\cM_{R}$ which can be interpreted as the
different twisted cohomology groups of $\cA//\cG$. But we will not go into
this further here.

\

\

Finally, we point out that, as explained by Burghelea and Haller, in the finite dimensional situation the twisted cohomology groups can be computed by a variant of Witten's complex, and described using differential forms satisfying Witten's deformation of the Hodge-Laplace equation.  It could well be that the structure we are discussing in this article arise this way in Quantum Field Theory, and are perhaps well-known there.

\section{Interaction between gauge theory and calibrated geometry}
\subsection{$G_{2}$-instantons and associative submanifolds}
So far we have ignored the fundamental problem of compactness  of moduli spaces. In this section we attempt to explore this and to get a glimpse of what modifications are  required to take account of the problem.

Suppose that $\sigma_{i}$ is a sequence of $G_{2}$-structures on fixed $7$-manifold $Y$ with limit $\sigma$. Suppose we have corresponding sequence of $G_{2}$-instantons on a fixed $SU(l)$ bundle. Then, according to Tian, after taking a subsequence the connections converge off a set $P$ of Haussdorf dimension at most  $3$ and each three-dimensional component satisfies the associative condition, in a generalised sense. The sequence of $4$-forms $\Tr(F_{A_{i}}^{2})$ allows us to attach a multiplicity to each $3$-dimensional component. Let us suppose that in fact $P$ is a connected smooth associative submanifold. Roughly speaking, we expect that if the multiplicity is $k$ then transverse to $P$ the connections  are modelled on instantons on $\bR^{4}$ with $c_{2}=k$,  and the behaviour as $i\rightarrow \infty$ mimics the familiar bubbling of instantons over $4$-manifolds.  We expect also that  the singularity at $P$ in the limiting connection is removable, so that limiting connection extends to a smooth $G_{2}$-instanton over $Y$\cite{kn:T2}. In this section we discuss the converse question. Let $\sigma$ be a $G_{2}$ structure on $Y$ and $A$ be  a  $G_{2}$-instanton  on a bundle $E$ over $Y$. Let $k$ be a positive integer and $P$ be an associative submanifold in $Y$. When does the triple $(A,P,k)$ appear as the limit of smooth $G_{2}$-instantons with respect to a  sequence of deformations
$\sigma_{i}$ of $\sigma$?

\

In this subsection we will explain
 that there is a natural candidate criterion for \lq\lq bubbling'' question. In particular when $l=2$ and $k=1$ we will argue that this occurs if
for some spin structure on $P$ the coupled Dirac operator on $E\vert_{P}$, defined by the restriction of the connection $A$, has a nontrivial kernel.
Several authors have considered related problems, mostly emphasising the similar question involving Cayley submanifolds and $\Spin(7)$-instantons in $8$-dimensions. Lewis proved an existence theorem for $\Spin(7)$-instantons
using a gluing construction, choosing a Cayley submanifold with very special properties\cite{kn:CL}. Brendle considered the  general question of existence of Yang-Mills connections \cite{kn:Br2}, and $\Spin(7)$-instantons in particular
\cite{kn:Br}, but restricting attention to the case when (the analogue of) the integer $k$ above is $1$. The construction we want to explain here is due to Haydys \cite{kn:Hay},\cite{kn:Hay2}, and related to ideas of Pidstrigatch \cite{kn:VP} and Taubes \cite{kn:Taub2}.  We refer to the forthcoming paper of Haydys \cite{kn:Hay2} for a more complete account, and a discussion of various other interesting related matters.

To begin, suppose that $V$ is a quaternionic manifold, with a multiplication map
\begin{equation} \mu:  TV\times \bH\rightarrow TV.\end{equation}
Then there is an elliptic \lq\lq Fueter equation'' for maps
$f:\bR^{3}\rightarrow V$ which is 
\begin{equation}  I \frac{\partial f}{\partial y_{1}} + J \frac{\partial f}{\partial y_{2}} + K\frac{\partial f}{\partial y_{3}}=0. \end{equation}
In the case when $V$ is $\bH$ this is just the Dirac equation for a spinor field. (There is a similar equation for maps from $\bR^{4}$ to $V$, but we will emphasise the $3$-dimensional version.) Now suppose that there is an action of $SU(2)$ on $V$ permuting the $I,J,K$. More precisely, this means that $\mu$ in (41) is an $SU(2)$-equivariant map, for the induced action of $SU(2)$ on $TV$ and the standard action by automorphisms of the quaternions. Let $P$ be an oriented Riemannian $3$-manifold with a spin structure and $Fr\rightarrow Y$ the corresponding principle $SU(2)$  bundle.  Then we can form the associated bundle \begin{equation}
\underline{V}= Fr\times_{SU(2)} V. \end{equation}
For each point $y\in P$ there is an obvious way to make $\bR\oplus TP_{y}$ into an algebra $H_{y}$ , isomorphic to $\bH$ but not canonically so. Our hypotheses imply that there is a natural $H_{y}$ structure on the tangent bundle of the fibre $\underline{V}_{y}$. Thus there is a Fueter equation for sections of $\underline{V}$,
$$  \sum_{i} e_{i} \nabla_{i} s=0, $$
where $e_{i}$ is any orthonormal frame in $TP_{y}$ and $\nabla_{i}$ denotes the \lq\lq vertical part'' of the derivative of the section, defined using the horizontal subspace induced from the Levi-Civita connection. Slightly more generally still, suppose that $G$ is another Lie group and there is an action of $G\times SU(2)$ on $V$, where now $G$ {\it preserves} the quaternionic structure. Let $Q\rightarrow P$ be a $G$-bundle with connection, so we have a principle $G\times SU(2))$ bundle $Q\times_{P}Fr$ over $Y$. Then we can form an associated fibre bundle 
$$  \left(Q\times_{P}Fr\right)\times_{G\times SU(2)} V, $$
and, as Haydys observed \cite{kn:Hay}, there is an intrinsic Fueter equation for sections of this bundle over $P$.

With this background in place, return to consider an associative submanifold $P\subset Y$. To explain the basic idea, begin with the model case when $P,Y$ are flat, so we can take local co-ordinates $y_{1}, y_{2}, y_{3}, x_{1},x_{2},
x_{3}, x_{4}$ in which 
$$  \sigma= dx_{1}dx_{2}dx_{3} dx_{4} + \sum_{i,j,k {\rm cyclic}} dy_{i}dy_{j} \omega_{k},$$
where $\omega_{i}$ is a standard basis for $\Lambda^{+}\bR^{4}$. For $\epsilon>0$ let $\sigma_{\epsilon}$ be $\epsilon^{-2}$ times the pull-back of $\sigma$ under the map which multiplies the $x_{a}$ co-ordinates by $\epsilon$. Thus
$$\sigma_{\epsilon}= \epsilon^{2} dx_{1}dx_{2} dx_{3}dx_{4}+\sum dy_{i} \omega_{i}. $$
We have a limit $\sigma^{*}=\sum dy_{i}dy_{j} \omega_{k}$ which is {\it not} a positive form, but we can still consider the equation $F\wedge \sigma^{*}=0$, which we refer to as the $\sigma^{*}$-instanton equation.

Let $N\rightarrow P$ be the normal bundle of $P\subset Y$. The $G_{2}$-structure induces a bundle isomorphism $\Lambda^{+}N\rightarrow TP$, which one finds is covariant constant with respect to the standard induced connections on $TP, N$. Fix a spin structure on $P$. Then we get another complex vector bundle $U\rightarrow P$, with  a connection and structure group $SU(2)$, such that $N$ is canonically identified with the real part of $U\otimes_{\bC}S$.
Using the connection on $N$ we get a canonical $3$-form $\sigma^{*}$ on the total space of $N$. Let $\exp$ be the exponential map from a neighbourhood of the zero-section in $N$ to $Y$, let $\exp_{\epsilon}(\xi)= \exp(\epsilon \xi)$ and let $\sigma_{\epsilon}$ be the $4$-form $\epsilon^{-2}\exp_{\epsilon}^{*}(\sigma)$.
Then one can see that the limit of $\sigma_{\epsilon}$ as $\epsilon$ tends to zero is $\sigma^{*}$. We can define $\sigma^{*}$ instantons:  connections on bundles over the total space $N$, as above and it clearly reasonable to expect that these are the blow-up limits of sequences of connections developing a singularity along $P$.

Now we want to bring in a $G_{2}$-instanton connection $A$ on another $SU(l)$ bundle $E$ over $Y$. This defines a connection $A\vert_{P}$on the restriction $E\vert_{P}$. Let $N_{\infty}$ be the $4$-sphere bundle over $P$ obtained by adjoining a section at infinity to the $\bR^{4}$ bundle $N$. W define a $\sigma^{*}$-instanton on $N$ with asymptotic limit $A\vert_{P}$ to be given by a  connection $\bA$ on a bundle $\bE\rightarrow N_{\infty}$ which satisfies the $\sigma^{*}$-instanton equations on the dense subset $N\subset N_{\infty}$ and such that the restriction of $\bA$ to the infinity section is equivalent to $A\vert_{P}$. Note that this data defines an integer  Chern class, given by the restriction of $\bE$ to any $4$-sphere fibre. We expect that  a triple
$(A,P, k)$ can only occur as the limit of a sequence of $G_{2}$ instantons, for nearby structures, if there is a solution of the $\sigma^{*}$- instanton equation with Chern class $k$ and asymptotic limit $A\vert_{P}$.

Let $M=M_{k,l}$ be the moduli space of \lq\lq framed'' $SU(l)$ instantons of Chern class $k$ over $S^{4}=\bR^{4}\cup\{\infty\}$. By a framing we mean a trivialisation of the fibre over $\infty$. We write $\bR^{4}$ as the  real part of $S^{+}\otimes S^{-}$ where $S^{+},S^{-}$ are the spin spaces. Then there is a natural action of the group
  $SU(2)\times SU(2)\times SU(l)$ on $M$, where the two $SU(2)$ factors act on $S^{+},S^{-}$,  and $SU(l)$ acts on the framing at infinity.  There is a quaternionic structure on $M$ which is preserved by $SU(l)$ and the second copy of $SU(2)$ and permuted by the first copy of  $SU(2)$. So we are in the situation above, with $G=SU(2)\times SU(l)$. Given our $SU(l)$ bundle $\bE\vert_{P}$ and our $SU(2)$ bundle $U\rightarrow P$ we form a  $G$-bundle as the fibre product and then we get a bundle
$\ccM\rightarrow Y$ with fibre $M$ and a Fueter equation for sections of $\ccM$. 

\begin{thm} (Haydys \cite{kn:Hay2})
There is a one-to-one correspondence between solutions of the $\sigma^{*}$
instanton equation with Chern class $k$ and asymptotic limit $A\vert_{P}$ and sections of the bundle $\ccM$
over $P$ which satisfy the Fueter equation.
\end{thm}

This can be thought of a variant of the \lq\lq adiabatic limit'' for $\Spin(7)$-instantons over products discussed in \cite{kn:DT}. In the case when $k=1$ the Fueter equation appears in \cite{kn:Br} as a \lq\lq balancing condition''. It is natural to expect that the equation in \cite{kn:Br2} associated to general Yang-Mills solutions can be interpreted as the equation defining a harmonic section of the bundle $\ccM$, as in \cite{kn:Hong}. 

\

We use the connection on $N$ to split the tangent space of $N$ into horizontal
and vertical subspaces, isomorphic to $TP$ and $(U\otimes S)_{\bR}$ respectively.
With respect to this splitting, $F(A)\wedge \sigma^{*}$ has two components,
say
$F_{1}, F_{2}$ where $F_{1}$ takes values in the normal bundle of $P$ and $F_{2}$ in the tangent bundle (both lifted to $N$ and tensored with the bundle of Lie algebras).  The $\sigma^{*}$-instanton equation thus splits into two separate conditions.
 We show that there is a one-to-one correspondence between
 \begin{itemize} \item connections
 $\bA$ over $N_{\infty}$ isomorphic to $A$ over the section at infinity and
 satisfying
  $F_{2}=0$
 \item   smooth sections of $\ccM\rightarrow Y$.
 \end{itemize}
 Then we show that the further condition $F_{1}=0$
 is equivalent to the Fueter equation.
 
 The condition $F_{2}=0$ just asserts that
 the restriction of $\bA$ to each fibre of $N$ is an anti-self-dual connection.
 There is a tautological bundle over $ M\times S^{4}$ which is equivariant
 for an  action of $SU(2)\times SU(2)\times SU(l)$ and which has a
 fixed trivialisation over $M\times \{\infty\}$. On this bundle we have a
 standard connection which restricts tautologically to the $S^{4}$ slices
 and which is compatible with the trivialisation over $M\times \{\infty\}$.
 Using the group action we construct a bundle $\tE$ over the pull-back
 $\pi^{*}(\ccM)$ of $\ccM$ to $N_{\infty}$, which can also be viewed as 
 a bundle over $P$ with fibre $S^{4}\times M$. The connections
 on $TP, U$ and $E$ induce a natural connection $\bA_{0}$ on $\tE$.
 A section $s$ of $\ccM$ induces a section $\overline{s}$ of $\pi^{*}(\ccM)$  and we have a connection $A_{s}= \overline{s}^{*}(\bA_{0})$ on a bundle
 $\overline{s}^{*}(\tE)$ over $N_{\infty}$. The stated properties of the connection
 on the universal bundle imply that $A_{s}$ is isomorphic to $A$ over the section at infinity
 and satisfies the condition $F_{2}=0$,
 simply because the connection on the universal bundle is anti-self-dual on each $\bR^{4}$
 slice in $M\times \bR^{4}$. Conversely, it is a straightforward formal exercise
 to check that all such connections arise in this way.

 The remaining task is to match up the curvature condition $F_{1}=0$ with the Fueter equation. For this we should recall
some more of the theory of instantons over $\bR^{4}$. We can
regard $\bR\oplus \Lambda^{+}$ as an algebra $H$. Of course, as before, $H$
is isomorphic to the quaternions but we prefer not to fix such an isomorphism.
Then $H$ acts naturally by multiplication on itself and also acts on $\bR^{4}$.
Let $A$ be a finite energy instanton over $\bR^{4}$. We have a defomation
complex
  $\Omega^{0}\stackrel{d_{A}}{\rightarrow} \Omega^{1} \stackrel{d^{+}_{A}}{\rightarrow} \Omega^{+}$
  and an elliptic operator $d^{*}_{A}\oplus d^{+}_{A}: \Omega^{1}\rightarrow
  \Omega^{0}\oplus \Omega^{+}$. The tangent space of the moduli space
  $M$ at $A$ can be identified with the $L^{2}$ solutions $a$ of $(d^{*}_{A}\oplus
  d^{+}_{A}) a=0 $ \cite{kn:Taub1}. The crucial points are 
  \begin{itemize}
  \item  $d^{*}_{A}\oplus d^{+}_{A}$
  commutes with the $H$ action induced by the actions on $\bR\oplus \Lambda^{+}$
  and $\bR^{4}$. Thus we get an action of $H$ on the tangent space of $M$, which
  is just the quaternionic structure mentioned before.
  \item The component of the curvature of the connection $A$ on the universal
  bundle in $T^{*}M\otimes T^{*}\bR^{4}={\rm Hom}(TM, T^{*}\bR^{4})$ is the tautological
  map given by evaluating $a\in TM$ at a point in $\bR^{4}$. In particular
  this commutes with the action of $H$. 
  \end{itemize}
  Now work at a fixed point in $P$ and  fix an orthonormal basis $e_{i}$ for the tangent
  space of $P$ at this point. Identifying the fibre of $N$ at this point
  with $\bR^{4}$ we get a basis $\omega_{i}$ of $\Lambda^{+}$. Suppose we
  have a section $s$ of $\ccM\rightarrow P$. With respect to the given connections
  this has a covariant derivative, with three components  $a_{1}, a_{2},
  a_{3}\in TM $  corresponding to the tangent vectors $e_{i}$.  The Fueter
  equation is $\sum \omega_{i}(a_{i})=0$. By the second observation above
  this implies that for each point $x\in \bR^{4}$
  $$  \sum \omega_{i} (F_{a_{i}), \ }=0, $$
  where $F_{a_{i},\ }$ is the bundle-valued  $1$-form on $\bR^{4}$ obtained
  by pairing the curvature of the universal bundle with $a_{i}\in TM$. Unravelling
  the definitions one sees that the left hand side of this equation  is precisely the component $F_{1}$ of the curvature. Thus a solution of the Fueter equation
does yield a solution of the $\sigma_{0}$-instanton equation and the only remaining
thing is to see that there are no other solutions.

Suppose we have any solution $\bA$ of the equation 
$F_{2}=0$, isomorphic to $A$ over the
infinity section. Restriction to the fibres of $N$ defines a section $s$ of $\ccM$ and $\bA$ must agree with $A_{s}$ in the fibre direction. Thus
the only possibility is that $\bA=A_{s}+ \Phi$ where $\Phi$ is a bundle-valued
$1$-form vanishing in the fibre direction. If we fix a point in $P$ and a
basis $e_{i}$ as above then $\Phi$ has three components $\phi_{i}$ which
are sections of the adjoint bundle over $\bR^{4}$. If $a_{i}$ are the derivatives
of $s$, as above, then the curvature condition $F_{1}=0$
becomes
\begin{equation}   \omega_{i}( a_{i} + d_{A} \phi_{i})=0. \end{equation}
The hypothesis on the connection over the infinity section is equivalent to
the condition that $\vert \phi_{i}(x)\vert \rightarrow 0$ as $x\in \bR^{4}$
tends to infinity. Thus what we  need to show is that in this situation all the $\phi_{i}$
vanish. Now write $\Phi=\sum \phi_{i} \omega_{i}$, a bundle-valued self-dual
$2$-form over $\bR^{4}$. The equation (44) is equivalent
to 
$$   d^{*}_{A} \Phi= \sum \omega_{i}( a_{i}). $$
Since $d^{+}_{A} \omega_{i} (a_{i})=0$ we have the identity $d^{+}_{A} d^{*}_{A} \Phi=0$.
The Weitzenbock formula in this situation tells us that
$$  d^{+}_{A} d_{A}^{*}\Phi= \sum_{i} \left( d_{A}^{*}d_{A} \phi_{i}\right)\omega_{i}.
$$
(This uses the fact that $A$ is an anti-self-dual connection.) So we deduce
that $ d_{A}^{*}d_{A}\phi_{i}=0$ and then the maximum principle implies that
$\phi$ vanishes, since $\vert \phi_{i}\vert $ tends to zero at infinity. This completes
the proof of the theorem.

\

There is a standard map from $M$ to $\bR^{4}$ which takes a connection to
the centre of mass of its curvature density $\vert F\vert^{2}$ and the derivative
of this map is $H$-linear. In fact $M$ is a product $M'\times \bR^{4}$, where
$M'$ is the \lq\lq centred'' moduli space. It follows that there is a bundle map from $\ccM$
to $N$ which takes solutions of the Fueter equation for sections of $\ccM$
to sections of the corresponding equation for sections of $N$. The latter
is just the linear Dirac equation appearing in the theory of deformations
of the associative submanifold $P$. We assume $P$ is \lq\lq regular''
so this equation has no non-zero solution. This means that we can replace
the instanton moduli space $M$ by the centred moduli space $M'$  throughout the discussion above. Let us consider the case when
$l=2$ and $k=1$. Then, up to translation and dilation of $\bR^{4}$ and gauge
equivalence, there is just one  instanton which is the standard connection
on the negative spin bundle over $S^{4}$. The framed moduli space can naturally
be written as
  $$M'= S^{+}\setminus \{0\} /\pm 1, $$
  where $S^{+}$ is the spin space, and this is compatible with the quaternionic
  structure. (The quotient by $\pm 1$ comes from the centre of $SU(2)$.) Tracing
  through the definitions we find that 
  $$  \ccM'=\left( (S^{+}\otimes E)_{\bR}\setminus \underline{0}\right)/\pm
  1. $$
  where $S^{+}$ is the spin bundle over $P$ and $\underline{0}$ denotes the
  zero-section. Note  that the bundle $\ccM'$ does not depend on the choice of
  a spin structure on $Y$ but if we have a section of $\ccM'$ there is a
  unique
  choice of spin structure for which this lifts to a section of $(S^{+}\otimes
  E)_{\bR}$. Making this choice, we see that solutions of the Fueter equation
  correspond (up to $\pm 1$) to solutions of the Dirac equation for sections
  of $E$-valued spinors over $P$, using the Dirac operator coupled to $A$.
  Thus, in sum we expect that the condition that a pair $P,A$ appear as the
  limit of $G_{2}$-instantons, with multiplicity $k=1$, is that there is a nontrivial harmonic spinor
  for the restriction $A\vert_{P}$, for some choice of spin structure on
  $P$. (Note that we are assuming here that the harmonic spinor does not vanish,
  but this should be true generically for dimensional reasons. A case when
  the harmonic spinor vanishes somewhere would require further analysis.)

  \

We 
   could carry out the same discussion for complex curves in a Calabi-Yau
  manifold and Hermitian-Yang-Mills connections, or embed this in the situation above by taking the product with
  a circle. In this 6-dimensional case  there is an alternative, algebro-geometric, point of view, which leads to the same conclusion. Let $\Sigma$ be a smooth curve in
a Calabi-Yau $3$-fold $Z_{0}$. Suppose $\pi:\cZ\rightarrow \Delta$ is a deformation of $Z_{0}$, where $\Delta$ is a disc in $\bC$. Thus $\Sigma$ can be viewed as a curve in the central fibre $\pi^{-1}(0)\subset \cZ$. Blow up this curve to get $\tpi: \tcZ\rightarrow \Delta$. Then
$$\tpi^{-1}(0)=  V\cup_{D} \tilde{Z}_{0}, $$
where $\tilde{Z}_{0}$ is the blow-up of $Z_{0}$ along $\Sigma$, $D\subset \tilde{Z}_{0}$ is
the exceptional divisor, which is a $\cP^{1}$ bundle over $\Sigma$, and $V$
is a $\cP^{2}$ bundle over $\Sigma$ which contains a copy of $D$. If we have
a holomorphic bundle $\bE$ over $\tcZ$, restriction to fibres gives a family of bundles $E_{t}$ over the deformations $Z_{t}=\pi^{-1}(t)$ for non-zero $t$ but
if the restriction of $\bE$ to the $\bC\bP^{2}$ fibres in $V$ is non-trivial
this family will not extend to give a bundle over $Z_{0}$. One expects this to give the algebro-geometric description
of a family of Hermitian-Yang-Mills connections  developing a singularity
along $\Sigma$. The analogue of the connection $A$ in the discussion above
is furnished by the restriction of $\bE$ to $\tilde{Z_{0}}\subset \tpi^{-1}(0)$,
which we assume to be the lift of a bundle $E$ over $Z_{0}$. The algebro-geometric
analogue of the question we have discussed above is to ask: given $\Sigma,
E, k$ when is there a bundle over $\tpi^{-1}(0)$ which is isomorphic to the
pull-back of $E$ over $\tilde{Z_{0}}$ and has $c_{2}=k$ on the $\bC\bP^{2}$ fibres
in $V$. It is a straightforward algebraic geometry exercise to show that
when $k=1$ this occurs precisely when there is a non-vanishing holomorphic section of
$E\otimes K^{1/2}_{\Sigma}$ over $\Sigma$,  for some choice of spin structure $K^{1/2}_{\Sigma}$.

\

\subsection{Implications}

 The Dirac operator on $E$-valued spinors over $P$
is naturally a real operator and we expect to encounter a zero eigenvalue
in real codimension $1$.  Thus  it seems likely that a naive count
of $G_{2}$-instantons will {\it not} yield an invariant. What one would expect
is needed is a count which includes triples $(A,P, k)$ of a connection on a different
bundle and an associative submanifold $P$, thought of as having multiplicity
$k$. We should count these with some weight $W(A,P,k)$. For example with
bundles of rank $l=2$ and when $k=1$ we would need some way to determine
the weight $W$ so that in a generic $1$-parameter family it changes by $\pm
1$ when one eigenvalue of the Dirac operator crosses $0$ (and a prerequisite
for understanding the signs would be to develop a theory of orientations
and signs in the \lq\lq naive'' theory, which we have not discussed).  This
is much the same issue involved in regularising the \lq\lq dimension'' in the ordinary Floer theory, as we discussed in Section 4.1.  Formally $W$ could be given by the $\pm 1/2$ where the sign is determined by the parity of the \lq\lq number of
negative eigenvalues'' of the coupled Dirac operator on $P$.  It seems possible
that this can be done, using the theory of spectral flow.

 For larger
values of $k$ new issues arise, since the contribution from $P$ involves an
essentially non-linear problem. We should still expect to encounter solutions
in real codimension $1$. The reason is that there is a dilation action on
$M$ and hence on $\ccM$ which preserves the quaternionic structure. Thus
 a single solution of the Fueter equation generates a $1$-dimensional family,
 by dilation. The linearisation of the Fueter equation has Fredholm index
 $0$ but this dilation action implies that, in a family, we we expect to encounter solutions
 in real codimension $1$, just as for the linear Dirac operator. Given any
 section $s$ of $\ccM$, let $\hat(s)$ denote the vertical vector field defined
 by the infinitesimal dilation action, and let ${\cal D}(s)$ be the expression
 appearing in the Fueter equation, which is also a vertical vector field.
 Then we have a nonlinear eigenvalue equation
 $$    {\cal D}(s)= \lambda \hat{s}, $$
 for sections $s$, generalising the eigenvalue equation for the Dirac operator.
 (In \cite{kn:Haythesis}, Section 3.4 , Haydys develops a more general theory of these eigenvalue equations, in terms of a \lq\lq Swann bundle''.)
 Of course we have a solution of the Fueter equation just when there is a
 zero eigenvalue. So it seems that one needs an extension of the theory
 of spectral flow which would enable one to define the weight $ W(A,P,k)$
 by a regularisation of the \lq\lq number of negative eigenvalues'' for this
 nonlinear problem.   If one seeks, more ambitiously, to construct
a Floer theory in $7$-dimensions then it  seems likely that one would have to assign a Floer group (or, perhaps better, chain complex) to $(A,P,k)$,
giving the contribution to the overall Floer homology (and with Euler characteristic
$W(A,P,k)$). This may be related to recent  work of Hohloch, Noetal and Salamon \cite{kn:DSal2}. We could think of a \lq\lq completed'' space of connections
with a point having a neighborhood modelled on the product of a Hilbert space
with a cone over a space ${\cal
L}$, where   $\cL$ is the space of sections of $\ccM$, modulo dilation. In
a finite-dimensional analogue the contribution of this point in a Morse theory
description of the homology of the total space will involve the theory of
the Conley index; the homology of $\cL$ and  various subsets. It is possible
that there is a \lq\lq Floer-analogue'' of this which can be formulated in terms of
the solutions of the eigenvalue equation, and \lq\lq flow lines'' between
them.

\

Even leaving aside this complication in $7$-dimensions, of bubbling along
an associative submanifold, it seems likely that the naive count of solutions
in $6$ dimensions is not the right thing to consider for the purposes of
developing a $6+1$ dimensional theory, including the conjectural holomorphic
bundle $\cE$ we have discussed. Consider a nonsingular  $G_{2}$
instanton over a topological cylinder $\bR\times Z$. Restriction  to slices $\{
t\} \times Z$ gives a family of connections over $Z$. One can imagine a case
when as $t\rightarrow\infty$ these connections develop a singularity along
$\Sigma\subset Z$. This phenomenon suggests that one would have to take account
of pairs $(E,\Sigma)$  in $Z$ in order to have the correct gluing identities.
This fits in with the fact that the Thomas' algebro-geometric approach to
the holomorphic Casson invariant includes contributions from sheaves, not
just bundles.  There seems to be a lot of  scope for work relating
the algebro-geometric and differential-geometric points of view.

\subsection{Codimension-3 theories and monopoles}

 This subsection is rather more speculative. Given a noncompact Riemannian $3$-manifold $B$  with an end of a suitable
 kind one can study the Bogomolny monopole equation $F(A)= *d_{A} \Phi$ for
 a connection $A$ and section $\Phi$ of the adjoint bundle. We  will just consider
 the case of structure group $SU(2)$. Solutions correspond
 to translation-invariant instantons on $\bR\times B$. One also imposes 
 asymptotic conditions on the ends of $B$, the most important being that $\vert \Phi\vert\rightarrow r^{-1}$ at infinity, where $r>0$ is fixed. The most familiar case is when $B=\bR^{3}$ and
 then one can reduce to the case when $r=1$ by scaling the metric. But in
 general $r$ will be a genuine parameter and it will not be possible to identify
 solutions for different values of $r$. The references \cite{kn:Braam}, \cite{kn:Floer1}, \cite{kn:JT}, \cite{kn:AH}, and many others, give more details about
this monopole theory.
 
 Now suppose we have a {\it noncompact} $G_{2}$-manifold $Y$. We can study
 the analogous equation (25) on $Y$ with the asymptotic condition $\vert u\vert\rightarrow
 r^{-1}$ at infinity. Let us imagine that, for manifolds $Y$ with an end of
 a suitable kind,  we can find a set-up which leads to a Fredholm problem  and to  invariants, which would be numbers in the case of index zero.  Then we could
study the behaviour of solutions as the parameter $r$ varies, in particular
as $r\rightarrow 0$. We  can find plausible models for this based on compact co-associative
submanifolds  $Q\subset Y$, in much the same way as we modelled the blow-up  behaviour
around associative submanifolds.  For simplicity consider first the flat
case, so we  have standard co-ordinates $x_{a}$ on $Q$ and $y_{i}$ normal to $Q$. The equations (25) can be written very schematically as
\begin{equation} \nabla_{\uy} u = F_{\uy \uy} + F_{\ux \ux}\ \ \  \nabla_{\ux} u = F_{\ux\uy}. \end{equation}
If we change variables, replacing $y_{i}$ by $r y_{i}$ and $u$ by $r^{-1} u$, then take the limit as $r$ tends to zero, we get a limiting equations
\begin{equation} \nabla_{\uy} u = F_{\uy\uy}\ \ \ \  \nabla_{\ux} u = F_{\ux\uy}.\end{equation}
The first of these equations is the  Bogomolony monopole equation on each
$\bR^{3}$ slice normal to $Q$ and the second equation is the Fueter equation
for the resulting map from $Q$ to the moduli space of monopoles on $\bR^{3}$.
To say this more systematically and generally, recall that for each integer $k>0$ we
have a moduli space $\Mon_{k}$ of \lq\lq centred'' monopoles of charge $k$ on $\bR^{3}$. (The charge is a topological invariant given by the degree of $\phi$ over a large sphere, in any trivialisation of the bundle.) This moduli
space is  a hyperkahler manifold of real dimension $4(k-1)$ and the rotations
of $\bR^{3}$ act on $\Mon_{k}$, permuting the complex structures. Given any
compact
coassociative submanifold $Q\subset Y$ we can form a bundle $\uMon_{k}\rightarrow
Q$ with fibre $\Mon_{k}$ much as before and there is a Fueter equation for
sections of $\uMon_{k}$. We interpret solutions of this Fueter equation as
possible asymptotic limits for sequences of solutions of (25) over $Y$ as the parameter $r$ tends to $0$. These solutions are
localised around $Q$ and should be very close to  reducible away from $Q$, with the structure
group reduced to $S^{1}\subset SU(2)$ by the \lq\lq Higgs field'' $u$ (which
would be almost covariant constant away from $Q$). 

We could make exactly the same discussion for a noncompact Calabi-Yau manifold
$Z^{6}$, and we will now switch our focus to this case as it is  simpler. We suppose the elementary topology of the set-up allows us to associate
a class in $H_{3}(Z)$ to our problem, analogous to the monopole charge in
$3$-dimensions. Then for each class $\kappa\in H_{3}(Z)$ we expect to have
a numerical invariant $n_{\kappa}$, counting solutions of (22). Our model for the asymptotic behaviour involves
a special Lagrangian submanifold $P\subset
Z$. For each $k$ we form a monopole bundle over $P$ and we have a Fueter
equation for sections. This has index $0$ and we imagine we can define a
number $w(k,P)$ by counting the solutions. Then we could hope to express
the number $n_{\kappa}$ as some kind of count of  special Lagrangian submanifolds $P$, weighted
by the numbers $w(k,P)$. The simplest guess is a formula of the shape
\begin{equation}  n_{\kappa}= \sum_{\kappa=\sum k_{i} [P_{i}]}\prod_{i} w(k_{i},P_{i}).\end{equation}

This  may be a rather crude approximation to the truth of the matter,
because we have not discussed what happens when two special Lagrangians intersect.
But in any case we could hope that there is some way of computing $n_{\kappa}$
from data localised around special Lagrangian submanifolds.

This picture, if it stands up to closer scrutiny, is rather similar to Taubes'
relation between the Seiberg-Witten  and Gromov invariants of a symplectic
four-manifold $W$\cite{kn:Taub3}. The moduli space of vortices on $\bR^{2}$ would play the role
in that case that the moduli space of monopoles does above. For a given \lq\lq
charge'' $k$ the moduli space of vortices is just the kth. symmetric product
of $\bR^{2}$ \cite{kn:JT}.   When the ambient
space is actually a complex  surface  we arrive in the realm of ordinary algebraic geometry. Given a curve $\Sigma \subset W$ a  section of the appropriate
\lq\lq vortex bundle''
corresponds to an infinitesimal deformation of the kth. order formal neighbourhood
of $\Sigma$, as a subscheme of $W$.

All of this discussion assumes that we can indeed find the correct set-up
to define numerical invariants $n_{\kappa}$.  The point we want to emphasise
is that, if this can be done, one might hope that these are {\it easier}
to define than the counts of special Lagrangians. Thus, by analogy, the Seiberg-Witten
invariants of a $4$-manifold are much easier to define than the Gromov invariants.
Then we could take (47) as the definition, or a guide to the definition, of
an invariant counting special Lagrangian submanifolds.


\end{document}